\documentclass[reqno]{amsart}
\usepackage{amsmath, amssymb ,amsthm, amsfonts, amsgen}
\usepackage{color}
\usepackage[final]{graphicx}
\usepackage{tikz}
\usepackage{bbm}

\numberwithin{equation}{section} 
\def\e{{\varepsilon}}
\def\O{{\Omega}}

\definecolor{vg}{rgb}{0.0, 0.40, 0.15}
\def\Rb{{\mathbb R}}

\def\L{{\mathcal {L}}}

\newtheorem{Theorem}{Theorem}[section]
\newtheorem{Lemma}[Theorem]{Lemma}
\newtheorem{Proposition}[Theorem]{Proposition}

\newtheorem{Remark}[Theorem]{Remark}
\newtheorem{Definition}[Theorem]{Definition}
\newcommand{{\rr}}{{\mathbb R}}

\newenvironment{@abssec}[1]{%
	\if@twocolumn
	\section*{#1}%
	\else
	\vspace{.05in}\footnotesize
	\parindent .2in
	{\upshape\bfseries #1. }\ignorespaces
	\fi}
{\if@twocolumn\else\par\vspace{.1in}\fi}
\newcommand\AMSname{AMS subject classifications}

\title[Integral representation for a relaxed optimal design problem in $BH(\O;\Rb^d)$]{Integral representation for a relaxed optimal design problem for 
non-simple grade two materials }
\author[A.~C.~Barroso]{Ana Cristina Barroso}
\address[A.~C.~Barroso]{Departamento de Matem\'atica and CMAFcIO, 
Faculdade de Ci\^encias da Universidade de Lisboa,
Campo Grande, Edif\' \i cio C6, Piso 1,
1749-016 Lisboa, Portugal}
\email{acbarroso@ciencias.ulisboa.pt}
\author[E.~Zappale]{Elvira Zappale}
\address[E.~Zappale]{Dipartimento di Scienze di Base ed Applicate per l'Ingegneria, Sapienza - Universit\`{a} di Roma, Via Antonio Scarpa, 16, 00161 Roma, Italy and
CIMA, Universidade de \'Evora, Portugal}
\email{elvira.zappale@uniroma1.it}

\makeindex

\begin{document}

\begin{abstract}
A measure representation result for a functional modelling optimal design problems for plastic deformations, under linear growth conditions, is obtained. 

Departing from an energy with a bulk term depending on the deformation gradient and its derivatives, as well as a perimeter term, the functional in question corresponds to the relaxation of this energy with respect to a pair $(\chi,u)$, where $\chi$ is the characteristic function of a set of finite perimeter and $u$ is a function of bounded hessian.


\par
\noindent2020 \AMSname: 49J45, 49Q20, 26B25
\end{abstract}

\maketitle

\tableofcontents

\section{Introduction}

The mathematical analysis of problems in materials science is based, in many instances, on the study of variational models involving higher order derivatives of the underlying fields. 

For example, energies of the type
\begin{equation}\label{genener}
\int_{\Omega}f_0(x, u(x), \nabla u(x), \nabla^2 u(x)) \, dx,
\end{equation}
set in the space of functions of bounded hessian, 
have been used to model certain problems in plasticity, see \cite{D, D1, D4, DM1, DM2}, to detect the behaviour of martensitic materials, see \cite{BL, BL2}, or in imaging problems as in, for example, \cite{BKP},
or, more generally, in the study of so called non-simple grade two materials as in \cite{Tou, Tou2} (see also \cite{CI2}). 
In all of these contexts, the search for equilibria leads naturally to the study of lower semi-continuity and relaxation properties for \eqref{genener},
where $f_0$ is the volume energy density and $u\equiv(u_1,\dots u_d)$ belongs to the space $BH(\Omega;\mathbb R^d)$ 
composed of integrable vector-valued functions for which all components $D^2_{jk}u_i$, 
$j,k = 1, \ldots,N$, $i=1,\dots, d$, of the hessian tensor $D^2 u := D(\nabla u)$ are bounded Radon measures and $\nabla^2 u$ stands for the absolutely continuous part, with respect to the Lebesgue measure, of the distributional second order derivative $D^2u$. We refer to \cite{FLMP, Mey, T, D}, among a wide literature. 

A further application of models involving measure higher order derivatives arises in connection with the field of structured deformations. The aim of this theory is to describe, at the macroscopic level, the effects of both smooth and non-smooth geometrical changes that a body may undergo, and which occur at different sub-macroscopic levels. The first-order theory of \cite{DPO} was extended in \cite{OP} to the second-order case, in order to describe curvature and bending effects associated with jumps in gradients. Relevant spaces in which to set the problem are $BH$, $SBH$ and $SBV^2$. Lower semi-continuity, relaxation and integral representation results have been obtained in the references \cite{BMMO, FHP, H}, to name but a few. 

Also motivated by second-order structured deformations, lower semi-continuity for 
\eqref{genener} assuming that $f_0$ is $2$-quasiconvex (cf. Definition \ref{2qcx}),
and satisfies a growth condition of order $p > 1$, in the last variable,
was proved in \cite{FLP2} for $u \in BH^p(\Omega;\mathbb R^d)$, meaning that $\nabla^2 u \in L^p$.

Regarding surface energies, the study of lower semi-continuity properties for 
$$ u \in SBH(\Omega;\mathbb R^d) \mapsto \int_{\Omega \cap S_{\nabla u}}
\varphi(x,(\nabla u)^+(x),(\nabla u)^-(x),\nu_{\nabla u}(x)) \, d\mathcal H^{N-1}(x),$$
was undertaken in \cite{SZ2}, whereas in \cite{BL} a derivation for a thin-film variational principle,  well-suited to apply finite element methods, and in which interfacial energy is modelled by the total variation of the deformation gradient $\nabla u$, is presented.

The strategy presented in Section \ref{jt}, based on the global method for relaxation in \cite{BFM}, allows us, as in \cite{FHP, H}, to handle a generalisation of the energy in \eqref{genener}, suited for the modelling of two-phase,
or multi-phase, \color{black} materials.
We comment on this case in Remark \ref{remGM}, indeed the main focus of this paper is on the case pertaining to the modelling of the optimal design for a two-components material, where each component, described by its characteristic function 
$\chi: \Omega \to \{0,1\}$, may contain both the austenite and the martensite phases. To penalise the transition of the interface between these two phases, a bulk energy depending on the second gradient of the deformation $u$ is added, modelled via two different energy densities $W_0$ and $W_1$, according to the component in question
(for the case of a single component and $W_0=W_1$ convex  we refer to \cite{BJ, LDM}, for a similar problem in the setting of thin structures described by means of smooth deformations, and to \cite{BL, BL2}, where discontinuities are allowed). 
A penalisation of the interfaces \color{black}  between the components is also added in the form of total variation of the field $\chi$, leading to an energy of the form
\begin{equation}\label{enernew} 
\int_{\Omega}\big[\chi(x) W_1(\nabla^2 u(x)) 
+ (1 - \chi(x))W_0(\nabla^2 u(x)) + W_2(x,\chi(x),u(x), \nabla u(x))\big] \, dx
+ |D \chi|(\Omega),
\end{equation}
see also \cite{RJ, S, SZ0, GZ1, GZ2} for related models in the smooth deformation setting, that is, when the continuous energy densities $W_0$ and $W_1$ satisfy super-linear growth conditions.
In our case, we assume that the energy densities have linear growth which, together with the lack of reflexivity of the space $L^1$,  naturally sets the equilibrium problem in the space $BH(\Omega;\mathbb R^d)$
(see, for example,  \cite{KS1, KS2, KS3, D, D1, D4}), obtaining, in particular, a more versatile limiting model that comprises also \cite{BL, BL2}.

As, under standard assumptions, the term involving $W_2$ in the above expression is a continuous perturbation,
it does not affect the analysis presented in this paper so from now on we will neglect it, stressing that, from Theorem \ref{main} below, the relaxation problem we consider, using the full expression for the energy in \eqref{enernew}, would lead to an integral representation comprised of the terms
\begin{align}\label{fing1}
&\int_A Q^2f(\chi(x), \nabla^2 u(x)) \, dx 
+  \int_A W_2(x,\chi(x),u(x), \nabla u(x)) \, dx \nonumber\\
& +	\int_{A \cap (S_\chi \cup S_{\nabla u})} 
g_1(x, \chi^+(x),\chi^-(x),u(x), (\nabla u)^+(x),(\nabla u)^-(x),\nu(x))\,
d\mathcal H^{N-1}(x) \\
& + \int_A(Q^2f)^{\infty}(\chi(x), \frac{d D^c(\nabla u) }{d |D^c(\nabla u)|}(x)) \, d|D^c(\nabla  u)|(x).\nonumber
\end{align}
Here, and in the sequel, to simplify the notation, we let 
$f: \{0,1\} \times \mathbb R^{d\times N \times N}_s \to [0,+\infty)$
be defined as
\begin{equation}\label{densityint}
	f\left( q,\xi\right)  :=q  W_1(\xi)+ (1-q)W_0(\xi), 
\end{equation}
$Q^2f$ is the  $2$-quasiconvex envelope of $f$ in the second variable
and, for a fixed $q \in\{0,1\}$, we recall that the recession function of $f$, in its second argument, is given by
\begin{equation}\label{recS}
	f^{\infty}(q,\xi) := \limsup_{t\rightarrow+\infty}\frac{f(q,t\xi)}{t}.
\end{equation}

Given these considerations, we are lead to the minimisation problem 
\begin{equation}\label{orpb} 
	\displaystyle \min_{(\chi,u)} \int_{\Omega}\big[\chi(x) W_1(\nabla^2 u(x)) 
	+ (1 - \chi(x))W_0(\nabla^2 u(x))\big]  dx + |D \chi|(\Omega), 
\end{equation}
where, usually, the volume fraction of each of the phases is prescribed,
leading to a constraint of the form
$\displaystyle \frac{1}{{\mathcal L}^N(\Omega)}\int_\Omega \chi(x) \, dx= \theta,  \;\theta \in (0,1)$. 
It is sometimes convenient to replace this constraint 
by inserting, instead, a Lagrange multiplier in the modelling functional
which, in our optimal design context, becomes \begin{equation}\label{functnoper}F(\chi,u;\Omega) : = 
\displaystyle \int_{\Omega}\left[\chi(x) W_1(\nabla^2 u(x)) 
+ (1 - \chi(x))W_0(\nabla^2 u(x))\right] \, dx +   \int_\Omega k \chi(x) dx 
+|D \chi|(\Omega).
\end{equation}

Although compactness holds in $BH(\Omega;\mathbb R^d)$ for functionals of the form \eqref{functnoper}, when extended as $+\infty$ in 
$BH(\Omega;\mathbb R^d) \setminus W^{2,1}(\Omega;\mathbb R^d)$, even when the last term is neglected, \color{black} the problem
of minimising 
\begin{equation*} 
	\displaystyle \int_{\Omega}\left[\chi(x) W_1(\nabla^2 u(x)) 
	+ (1 - \chi(x))W_0(\nabla^2 u(x))\right] \, dx +   \int_\Omega k \chi(x) dx, \, \color{black}
\end{equation*}
with respect to $(\chi,u)$, adding suitable forces and/or boundary conditions, is known to be ill-posed, in the sense that
minimising sequences $\{\chi_n\} \subset L^\infty(\Omega;\{0,1\})$ tend to highly oscillate and develop microstructure, so that in the limit we may no longer obtain a characteristic function. 
As in \cite{ABper} and \cite{KL}, this phenomenon is avoided by considering  the full expression given in \eqref{functnoper}, i.e. by  adding the perimeter penalisation along the interface between the two zones 
$\{\chi = 0\}$ and $\{\chi = 1\}$ (we refer to \cite{CZ} for the analogous analysis performed in $BV$, and to \cite{CZ0, BZ, BZ2} for the Sobolev settings, also in the presence of a gap in the growth exponents). 

We also mention the work of Arroyo-Rabasa in \cite{AR} where,  in the same spirit of \cite{ABper}, a vast class of optimal design problems with a perimeter penalisation is studied. These are based on energies of the form

$$(u,A) \mapsto \int_{\Omega}fu \, dx - 
\int_{\Omega \cap A}\sigma_1 \mathbb B u \cdot \mathbb B u \, dx -
\int_{\Omega \setminus A}\sigma_2 \mathbb B u \cdot \mathbb B u \, dx +
{\rm Per}(A;\overline \Omega),$$
where $\sigma_1, \sigma_2$ are two symmetric tensors and, adopting the notation of \cite{Ra}, $\mathbb B$   belongs to a class of operators of gradient form which includes both scalar and vector-valued
gradients, symmetrised gradients and higher-order gradients. The results in 
\cite{AR}, however, do not cover our present case, as the problem is set only in the space
of functions $u \in L^2$ such that $\mathbb  Bu \in L^2$ and the functional is {\it convex}. We refer to Proposition \ref{casep} for an existence result in the 
{\it non-convex}, case within the Sobolev framework.

With the abuse of notation of denoting, simply by $W_1$, the function $W_1+ k$ in \eqref{functnoper}, 
\color{black} 
our aim in this paper is to study the energy functional given by
\begin{equation}\label{Fint}
F(\chi,u;\Omega) : =\begin{cases}
\displaystyle \int_{\Omega}\big[\chi(x) W_1(\nabla^2  u(x)) 
+ (1 - \chi(x))W_0(\nabla^2 u(x))\big] \, dx + |D\chi|(\Omega), & 
\hbox{ if } \chi \in BV(\Omega;\{0,1\}) \\
& \mbox{ and } u \in  W^{2,1}(\Omega;\mathbb R^d), \\
+\infty, &\hbox{otherwise,}\end{cases}
\end{equation}
where $u \in BH(\Omega;\mathbb R^d)$, $\chi \in BV(\Omega;\{0,1\})$ and the densities $W_i$, $i = 0,1$, are continuous functions satisfying the following linear growth conditions from above and below,
 \begin{equation}\label{growthint}
 	\exists \, \alpha, \beta > 0 \text{ such that }
 	\alpha |\xi| \leq W_i(\xi) \leq \beta (1 + |\xi|), \; \; \forall 
 	\xi \in \mathbb R^{d\times N \times N}_s.
 \end{equation}
Furthermore, we point out that no convexity assumptions are placed on 
$W_i$, $i = 0,1$.

It is an immediate consequence of \eqref{growthint} that 
\begin{equation}\label{G}
	|f(q_1,\xi) - f(q_2,\xi)| \leq \beta \, |q_1 - q_2|(1 + |\xi|), \; \;
	\forall q_1, q_2 \in \{0,1\},
	\forall \xi \in \Rb^{d \times N \times N}_s,
\end{equation}
from which it follows that 
\begin{equation*}
	|f^{\infty}(q_1,\xi) - f^{\infty}(q_2,\xi)| \leq \beta \, |q_1 - q_2| \,|\xi|, \; \;
	\forall q_1, q_2 \in \{0,1\},
	\forall \xi \in \Rb^{d\times N \times N}_s.
\end{equation*}

Since we place no convexity assumptions on $W_i$, we consider
the relaxed localised functionals arising from the energy \eqref{Fint}, 
defined, for an open subset $A \subset \Omega$, $\chi \in BV(\Omega;\{0,1\})$ and 
$u \in BH(\Omega;\mathbb R^d)$, by
\begin{align}\label{calFint}
\mathcal{F}\left(\chi,u;A\right)  
:=\inf\bigg\{ \liminf_{n\rightarrow +\infty} F(\chi_n,u_n;A) \, : & \, 
\{\chi_n\} \subset BV(A;\{0,1\}), \{u_n\} \subset W^{2,1}(A;\mathbb{R}^{d}),
\nonumber\\
&  \chi_{n}\to\chi\text{ in } L^1(A;\{0,1\}), 	
u_{n}\to u\text{ in } W^{1,1}(A;\mathbb{R}^{d})
\bigg\}.
\end{align}
 
Due to the expression of \eqref{Fint}, and to the fact that 
$\chi_n \overset{\ast}{\rightharpoonup}\chi$ in $BV$ if and only if $\{\chi_n\}$ is uniformly bounded in $BV$ and $\chi_n \to \chi$ in $L^1$, it is equivalent to take
$\chi_n \overset{\ast}{\rightharpoonup}\chi$ in $BV$ or $\chi_n \to \chi$ in $L^1$
in the definition of the functional \eqref{calFint} obtaining the same infimum regardless of the considered convergence.
 
We prove in Proposition~\ref{traceRm} that $\mathcal{F}\left(\chi,u;\cdot\right)$ is the restriction to the open subsets of $\Omega$ of a Radon measure, the main result of our paper concerns the characterisation of this measure.

\begin{Theorem}\label{main}
Let $f:\{0,1\} \times \mathbb R^{d\times N \times N}_s\to [0, + \infty)$ be a continuous function as in \eqref{densityint}, where $W_0$ and $W_1$ satisfy \eqref{growthint}, and consider 
$F:BV(\Omega;\{0,1\})\times BH(\Omega;\mathbb R^d)\times \mathcal O(\Omega) \to [0,+\infty]$ 
defined in \eqref{Fint}. Then, the functional in \eqref{calFint} admits the following representation
\begin{align*}
\mathcal{F}\left(  \chi,u;A\right)&=\int_A Q^2f(\chi(x), \nabla^2 u(x)) \, dx  \nonumber \\
& +	\int_{A \cap (S_\chi \cup S_{\nabla u})} 
g(\chi^+(x),\chi^-(x),(\nabla u)^+(x),(\nabla u)^-(x),\nu(x))\,
d\mathcal H^{N-1}(x) \nonumber  \\
& + \int_A(Q^2f)^{\infty}(\chi(x), \frac{d D^c(\nabla u) }{d |D^c(\nabla u)|}(x)) \, d|D^c(\nabla  u)|(x), 
\end{align*}
where $Q^2f$ is the  $2$-quasiconvex envelope of $f$ in the second variable and 
$(Q^2f)^{\infty}$ is its strong recession function (cf. Subsection~\ref{defqcx} and \eqref{recS}, respectively). The relaxed surface energy density is given by
\begin{equation}\label{1gdef}g(a,b,C,D,\nu) := \limsup_{\varepsilon \to 0^+}
\frac{m(\chi_{a,b,\nu},u_{C,D,\nu};Q_\nu(0,\e))}
{\e^{N-1}}
\end{equation}\color{black}
where $Q_{\nu}(0,\e)$ stands for an open cube with centre $0$, side-length $\e$ and two of its faces perpendicular to the unit vector $\nu$, 
$m: BV(\Omega;\{0,1\})\times BH(\Omega;\mathbb R^d)\times \mathcal O(\Omega) \to [0,+\infty]$ is defined by
\begin{align*}
m(\chi,u;A) := \inf\left\{\mathcal F (\theta,v;A) : 
\theta \in  BV(\O;\{0,1\}), v \in BH(\O;\mathbb R^d), \theta\equiv \chi 
\text{ on }\partial A, {\rm supp}(u-v) \Subset A \right\},
\end{align*}
and, for $(a,b,C,D,\nu)  \in \{0,1\} \times \{0,1\} \times \mathbb{R}^{d\times N} \times\mathbb{R}^{d\times N} \times \mathbb S^{N-1}$, with  $C-D= P \otimes \nu$, for some $P \in \mathbb R^d$, the functions $\chi_{a,b,\nu}$ and $u_{C,D,\nu}$ are defined as
\begin{equation*}
\chi_{a,b,\nu}(y) := \begin{cases}
a, & {\rm if } \; y \cdot \nu \geq 0\\
b, & {\rm if } \; y \cdot \nu < 0
\end{cases}
\; \; \; {\rm and } \; \; \;
u_{C,D,\nu}(y) := \begin{cases}
C \cdot y, & {\rm if } \; y \cdot \nu \geq 0\\
D \cdot y, & {\rm if } \; y \cdot \nu < 0.
\end{cases}
\end{equation*}
\end{Theorem} 

We point out that the condition $C-D= P \otimes \nu$, for some $P \in \mathbb R^d$, is required due to the symmetry of the Hessian and its rank-one property (see \cite[(2.2)]{FLP} and \cite{FLP2}). 

If the original bulk energy takes the more general form
$$\int_A f_1\color{black}(x,\chi(x),u(x),\nabla u(x),\nabla^2u(x)) \, dx,$$ 
which includes, in particular, the bulk energy in \eqref{enernew}, 
then the relaxed surface energy density in \eqref{fing1} may depend explicitly on $x$ and $u(x)$ and is given by 
$$g_1\color{black}(x_0,a,b,r,C,D,\nu) := \limsup_{\varepsilon \to 0^+}
\frac{m(\chi_{a,b,\nu}(\cdot - x_0),r+ u_{C,D,\nu}(\cdot - x_0);Q_\nu(x_0,\e))}
{\e^{N-1}},$$
where the cube $Q_{\nu}(x_0,\e)$ is now centred at $x_0$.
However, in the particular case of the energy in \eqref{Fint}, it turns out that the function $g$ is, in fact, independent of both $x$ and $u(x)$ since 
$$g_1(x_0,a,b,r,C,D,\nu) = g_1(0,a,b,0,C,D,\nu)$$
which, we denote by $g(a,b,C,D,\nu)$. More details may be found in Proposition \ref{mg} and Remark \ref{remGM}.

For the notation regarding the singular sets $S_\chi$, $S_{\nabla u}$ and the corresponding vectors
$\chi^+(x)$, $\chi^-(x)$, $\nu_\chi(x)$, $(\nabla u)^+(x)$, $(\nabla u)^-(x)$ and $\nu_{\nabla u}(x)$ we refer to Subsection~\ref{BV}. 
 
The above expression for the relaxed surface energy density, $g$, arises as an application of the global method for relaxation, see \cite{BFM, FHP}. 
However, as we will see in Proposition~\ref{Ktilde}, in the case where $f$ satisfies the additional hypothesis 
\begin{equation}\label{finfty}
\exists \, 0 < \gamma \leq 1, \exists \, C, L > 0 : \; t\,|\xi| > L \Rightarrow 
\left| f^{\infty}(q,\xi) - \frac{f(q,t\xi)}{t}\right| 
\leq C \frac{|\xi|^{1-\gamma}}{t^\gamma},
\end{equation}
for every $q \in \{0,1\}$ and every $\xi \in \mathbb R^{d \times N \times N}_s$
(which, as pointed out in \cite{FM}, can be stated equivalently as
\begin{equation}\label{finfty2}
\exists \, 0 < \gamma \leq 1, \exists \, C > 0 \mbox{ such that }
\left| f^{\infty}(q,\xi) - f(q,\xi)\right| 
\leq C \left(1 + |\xi|^{1-\gamma}\right),
\end{equation}
for every $q \in \{0,1\}$ and every $\xi \in \mathbb R^{d\times N \times N}_s$),
we can provide an alternative, sequential, characterisation of $g$, that is  independent of the Dirichlet functional \eqref{mbis}. Indeed, in this case, we can show that
$$g(a,b,C,D,\nu) = \widetilde{K}(a,b,C,D,\nu),$$
for every $(a,b,C,D,\nu)  \in \{0,1\} \times \{0,1\} \times \mathbb{R}^{d \times N} \times\mathbb{R}^{d \times N} \times \mathbb S^{N-1}$, such that $C-D=P\otimes \nu$, for some $P\in\mathbb R^d$, where
\begin{align*}
\widetilde{K}(a,b,C,D,\nu) & := \inf\bigg\{\liminf_{n \to + \infty}
\left[\displaystyle\int_{Q_{\nu}}f^{\infty}(\chi_n(x),\nabla^2 u_n(x)) \, dx
+ |D\chi_n|(Q_{\nu})\right] : \{\chi_n\} \subset BV\left(Q_{\nu};\{0,1\}\right),\\
& \hspace{1,1cm} \{u_n\} \subset W^{2,1}\left(Q_{\nu};\mathbb{R}^d\right),  
\chi_n \to \chi_{a,b,\nu}
\; {\rm in } \; L^1(Q_{\nu};\{0,1\}), u_n \to u_{C,D,\nu} \; {\rm in } \; W^{1,1}(Q_{\nu};\mathbb R^d) \bigg\}, 
\end{align*} 
and $\chi_{a,b,\nu}$, $u_{C,D,\nu}$ were defined in Theorem \ref{main}. In particular, $g$ is lower semi-continuous.

Furthermore, under suitable assumptions, we will show that $g$ is 
$BV\times SBH$-elliptic, see Proposition \ref{prop4.12}.

In the Appendix we present another, partial, more explicit, estimate for $g$ under hypothesis \eqref{finfty}. Here we prove that
$$g(a,b,C,D,\nu) \geq K(a,b,C,D,\nu)$$
where
\begin{equation}\label{Klb}
K(a,b,C,D,\nu):=\inf\left\{  
\displaystyle\int_{Q_{\nu}}\!\!\!
(Q^2f)^{\infty}(\chi(x), \nabla^2 u(x)) \, dx+|D\chi|(Q_{\nu}):\left(\chi,u\right)
\in\mathcal{A}(a,b,C,D,\nu)\right\}, 
\end{equation}
and, for $(a,b,C,D,\nu)  \in \{0,1\} \times \{0,1\}  \times \mathbb{R}^{d\times N} \times\mathbb{R}^{d\times N} \times \mathbb S^{N-1},$ such that $C-D= P \otimes\nu$, for some $P \in \mathbb R^d$, the set of admissible functions is
\begin{align*}
\mathcal{A}(a,b,C,D,\nu) &:=\bigg\{ \left(\chi,u\right)\in
BV\left(Q_{\nu};\{0,1\}\right)  \times 
W^{2,1}\left(Q_{\nu};\mathbb{R}^{d}\right): \chi= \chi_{a,b,\nu} \text{ on } \partial Q_\nu, u = u_{C,D,\nu}
\text{ on }\partial Q_\nu\bigg\}. 
\end{align*}

The fact that we are dealing with $BH$ fields, which allow only for jump discontinuities of their gradients, prevents us from obtaining the reverse inequality. In any case, Theorem \ref{main} generalises, to the optimal design scenario in the $BH$ setting, a similar integral representation result for \eqref{calFint} that was obtained in $BV$ in \cite{CZ}, when $W_0$ and $W_1$ depend on the gradient $\nabla u$. 
 
We remark that, in view of Theorem \ref{interpolBH} below, the integral representation, obtained in Theorem \ref{main}, also holds when the functional in \eqref{calFint} is replaced by 
\begin{align*}
 	\overline{\mathcal{F}}\left(\chi,u;A\right)  
 	:=\inf\bigg\{ \liminf_{n\rightarrow +\infty} F(\chi_n,u_n;A) \, : & \, 
 	\{\chi_n\} \subset BV(A;\{0,1\}), \{u_n\} \subset W^{2,1}(A;\mathbb{R}^{d}),
 	\nonumber\\
 	&  \chi_{n}\to\chi\text{ in } L^1(A;\{0,1\}), 	
 	u_{n}\to u\text{ in } L^{1}(A;\mathbb{R}^{d})
 	\bigg\},
\end{align*}
whenever the open subset $A \subset \Omega$ has Lipschitz boundary, 
$\chi \in BV(\Omega;\{0,1\})$ and $u \in BH(\Omega;\mathbb R^d)$. 

The case where the densities $W_i$, $i=0,1$, satisfy \eqref{growthint}, but the functional $\mathcal F$ in \eqref{calFint} is replaced by the the sequential lower semi-continuous envelope of  $F$ in \eqref{Fint}, with respect to the weak $W^{2,1}$ convergence, will be discussed in Proposition \ref{Frel=Frel**} below.
 
We conclude this introduction by observing that an analogue of Theorem \ref{main} in the case where the density functions $W_i$, $i=0,1$,  satisfy a growth condition of order $p>1$, i.e. when \eqref{growthint} is replaced by 
\begin{equation}
	\label{growthintp}
	\exists \, 0<\alpha \leq \beta  \hbox{ such that }
	\alpha |\xi|^p \leq W_i(\xi) \leq \beta(1+ |\xi|^p), \forall \xi \in \mathbb R^{d\times N \times N}_s,
\end{equation} can be easily deduced by existing relaxation results, see \cite[Theorems 1.1 and 1.3]{BFL}.

In fact, under assumption \eqref{growthintp}, the functional $\mathcal F(\chi,u;\Omega)$ in \eqref{calFint} coincides with
\begin{align*}
\mathcal{F}\left(\chi,u;\Omega\right) :=
\inf\bigg\{ \liminf_{n\rightarrow +\infty} F(\chi_n,u_n;\Omega) \, : & \,
\{\chi_n\} \subset BV(\Omega;\{0,1\}), \{u_n\} \subset W^{2,p}(\Omega;\mathbb{R}^{d}),\\
& \chi_{n}\to\chi\text{ in } L^1(\Omega;\{0,1\}), 
u_{n}\to u\text{ in } W^{1,1}(\Omega;\mathbb{R}^{d})\bigg\},
\end{align*}
which in turn, is finite only on 
$BV(\Omega;\{0,1\})\times W^{2,p}(\Omega;\mathbb R^d)$ and can be represented as
\begin{align*}
\mathcal F(\chi,u;\Omega)=
\int_\Omega Q^2f(\chi(x),\nabla^2 u(x)) \, dx + |D \chi|(\Omega).
\end{align*}
We refer to Proposition \ref{casep} for more details.
	
\medskip

The paper is organised as follows: in Section \ref{prelim} we set our notation and state some results on functions of bounded variation and on functions of bounded hessian that will be used in the sequel. Section \ref{auxres} contains some auxiliary results, namely some properties of the localised functional \eqref{calFint}, whereas Section \ref{mainsec} is devoted to the proof of the main theorem. In the Appendix we provide a lower bound for the relaxed surface energy density in terms of \eqref{Klb}.

\section{Preliminaries}\label{prelim}

In this section we fix the notation and quote some definitions and 
results that will be used in the sequel.

\subsection{Notation}\label{not}

Throughout the text $\Omega \subset \mathbb R^{N}$ will denote an open, bounded set
with Lipschitz boundary.

We will use the following notations:
\begin{itemize}
	\item ${\mathcal B}(\Omega)$, ${\mathcal O}(\Omega)$ and ${\mathcal O}_{\infty}(\Omega)$ represent the families of all Borel, open and open subsets of $\Omega$ with Lipschitz boundary, respectively;
\item $\mathcal M (\Omega)$  and $\mathcal M(\Omega;\mathbb R^h)$, $h \in \mathbb N$, $h > 1$, are the sets of finite real and vector-valued Radon measures on $\Omega$, respectively; 
	\item $\left |\mu \right |$ stands for the total variation of a measure  $\mu$ in $\mathcal M (\Omega)$ or in $\mathcal M(\Omega;\mathbb R^h)$;  
	\item $\mathcal L^{N}$ and $\mathcal H^{N-1}$ are the  $N$-dimensional Lebesgue measure 
	and the $\left(N-1\right)$-dimensional Hausdorff measure in $\mathbb R^N$, respectively;
	\item the symbol $d x$ will also be used to denote integration with respect to $\mathcal L^{N}$;
	\item the set of symmetric $N \times N$ matrices is denoted by $\Rb_s^{N\times N}$;
	\item if $u:\Omega \to \mathbb R$ is smooth,  its hessian matrix is an element of  $\Rb_s^{N\times N}$; 
	\item analogously, if  $u=(u^1,...,u^d):\Omega\to\Rb^d$ is smooth, then $D^2u$ is written as a vector in $\Rb^d$ whose $i$th component is the symmetric matrix 
	$\displaystyle \left(\frac{\partial^2 u^i}{\partial x_k\partial x_j}\right)_{1\leq k,j\leq N}$ in $\Rb_s^{N\times N}$, i.e. $D^2u \in (\Rb_s^{N\times N})^d$; however, we observe that in the framework of $BH$,  as in \cite{FLP} for instance,  $D^2u$ is often written as a third-order symmetric tensor, i.e. as an element of $\Rb_s^{d \times N\times N}$;
	\item $B(x, \e)$ is the open ball in $\Rb^N$ with centre $x$ and radius $\e$, 
	$Q(x,\e)$ is the open cube in $\Rb^N$ with two of its faces orthogonal
		 to the unit vector $e_N$, centre $x$ and side-length $\e$, whereas
	$Q_{\nu}(x,\e)$ stands for a cube with two of its faces orthogonal
	to the unit vector $\nu$; when $x=0$ and $\e = 1$, we write $Q_{\nu}$; if also $\nu=e_N$ we simply write $B$ and $Q$;
	\item $\mathbb S^{N-1} := \partial B$ is the unit sphere in $\Rb^N$;
	\item for every $h\in \mathbb N$, 
	$C_c^{\infty}(\Omega; \Rb^h)$ is the space of $\Rb^h$-valued smooth functions with compact support in $\Omega$;
	\item by $\displaystyle \lim_{k,n}$ we mean 
	$\displaystyle \lim_{k \to +\infty}\lim_{n \to +\infty}$;
	\item  $C$ represents a generic positive constant that may change from line to line.
\end{itemize}

\subsection{BV fields and sets of finite perimeter}\label{BV}

In the following we give some preliminary notions regarding functions of bounded variation and sets of finite perimeter. For a detailed treatment we refer to \cite{AFP}.

\smallskip

Given $u \in L^1(\Omega; \Rb^d)$ we let $\O _u$ be the set of Lebesgue points of $u$,
i.e., $x\in \O _u$ if there exists $\widetilde u(x)\in {\mathbb{R}}^d$ such
that 
\begin{equation*}
	\lim_{\varepsilon\to 0^+} \frac{1}{\e^N}\int_{B(x,\varepsilon)}
	|u(y)-\widetilde u(x)|\,dy=0,
\end{equation*}
$\widetilde u(x)$ is called the approximate limit of $u$ at $x$.
The Lebesgue discontinuity set $S_u$ of $u$ is defined as 
$S_u := \Omega \setminus \O _u$. It is known that ${\mathcal{L}}^{N}(S_u) = 0$ and the function 
$x \in \Omega \mapsto \widetilde u(x)$, which coincides with $u$ $\mathcal L ^N$- a.e.
in $\O _u$, is called the Lebesgue representative of $u$.

The jump set of the function $u$, denoted by $J_u$, is the set of
points $x\in \O \setminus \O _u$ for which there exist 
$a, \,b\in {\mathbb{R}}^d$ and a unit vector $\nu \in \mathbb S^{N-1}$, normal to $J_u$ at $x$, such that $a\neq b$ and 
\begin{equation*}  
	\lim_{\varepsilon \to 0^+} \frac {1}{\varepsilon^N} \int_{\{ y \in
		B(x,\varepsilon) : (y-x)\cdot\nu > 0 \}} | u(y) - a| \, dy = 0,
	\qquad
	\lim_{\varepsilon \to 0^+} \frac {1}{\varepsilon^N} \int_{\{ y \in
		B(x,\varepsilon) : (y-x)\cdot\nu < 0 \}} | u(y) - b| \, dy = 0.
\end{equation*}
The triple $(a,b,\nu)$ is uniquely determined by the conditions above,  
up to a permutation of $(a,b)$ and a change of sign of $\nu$,
and is denoted by $(u^+ (x),u^- (x),\nu_u (x)).$ The jump of $u$ at $x$ is defined by
$[u](x) : = u^+(x) - u^-(x).$

\smallskip

We recall that a function $u\in L^{1}(\Omega;{\mathbb{R}}^{d})$ is said to be of bounded variation, and we write $u\in BV(\Omega;{\mathbb{R}}^{d})$, if
all its first order distributional derivatives $D_{j}u_{i}$ belong to 
$\mathcal{M}(\Omega)$ for $1\leq i\leq d$ and $1\leq j\leq N$.

The matrix-valued measure whose entries are $D_{j}u_{i}$ is denoted by $Du$
and $|Du|$ stands for its total variation.
The space $BV(\O ; {\mathbb{R}}^d)$ is a Banach space when endowed with the norm 
\begin{equation*}
	\|u\|_{BV(\O ; {\mathbb{R}}^d)} = \|u\|_{L^1(\O ; {\mathbb{R}}^d)} + |Du|(\O )
\end{equation*}
and we observe that if $u\in BV(\Omega;\mathbb{R}^{d})$ then $u\mapsto|Du|(\Omega)$ is lower semi-continuous in $BV(\Omega;\mathbb{R}^{d})$ with respect to the
$L_{\mathrm{loc}}^{1}(\Omega;\mathbb{R}^{d})$ topology.

By the Lebesgue Decomposition Theorem, $Du$ can be split into the sum of two
mutually singular measures $D^{a}u$ and $D^{s}u$, the absolutely continuous
part and the singular part, respectively, of $Du$ with respect to the
Lebesgue measure $\mathcal{L}^N$. By $\nabla u$ we denote the 
Radon-Nikod\'{y}m derivative of $D^{a}u$ with respect to $\mathcal{L}^N$, so that we
can write 
\begin{equation*}
	Du= \nabla u \mathcal{L}^N \lfloor \O + D^{s}u.
\end{equation*}

If $u \in BV(\O )$ it is well known that $S_u$ is countably $(N-1)$-rectifiable, see \cite{AFP},  
and the following decomposition holds 
\begin{equation*}
	Du= \nabla u \mathcal{L}^N \lfloor \O + [u] \otimes \nu_u {\mathcal{H}}^{N-1}\lfloor S_u + D^cu,
\end{equation*}
\noindent where $D^cu$ is the Cantor part of the
measure $Du$.  When $D^cu = 0$, the function $u$ is said to be a special function of bounded variation, written
$u \in SBV(\Omega;\mathbb R^d)$. 

If $\Omega$ is an open and bounded set with Lipschitz boundary then the
outer unit normal to $\partial \Omega$ (denoted by $\nu$) exists ${\mathcal{H}}^{N-1}$-a.e. and the trace for functions in $BV(\Omega;{\mathbb{R}}^d)$ is
defined.

\begin{Theorem}
	(Approximate Differentiability in $BV$)\label{approxdiff} 
	If $u \in BV(\Omega; {\mathbb{R}}^d),$ then for $\mathcal{L}^N$-a.e. $x \in\Omega$ 
	\begin{equation*}
		\lim_{\varepsilon \rightarrow 0^+} \frac{1}{\varepsilon^{N+1} }
		\int_{Q(x, \varepsilon)} |u(y) - u(x) - \nabla u(x)\cdot(y-x)|\, dy  = 0.
	\end{equation*}
\end{Theorem}

\begin{Definition}
	\label{Setsoffiniteperimeter} Let $E$ be an $\mathcal{L}^{N}$- measurable
	subset of $\mathbb{R}^{N}$. For any open set $\Omega\subset\mathbb{R}^{N}$ the
	{\em perimeter} of $E$ in $\Omega$, denoted by $P(E;\Omega)$, is 
	given by
	\begin{equation*}
		P(E;\Omega):=\sup\left\{  \int_{E} \mathrm{div}\varphi(x) \,dx:
		\varphi\in C^{1}_{c}(\Omega;\mathbb{R}^{N}), \|\varphi\|_{L^{\infty}}%
		\leq1\right\}  .
	\end{equation*}
	We say that $E$ is a {\em set of finite perimeter} in $\Omega$ if $P(E;\Omega) <+
	\infty.$
\end{Definition}

Recalling that if $\mathcal{L}^{N}(E \cap\Omega)$ is finite, then $\chi_{E}
\in L^{1}(\Omega)$, by \cite[Proposition 3.6]{AFP}, it follows
that $E$ has finite perimeter in $\Omega$ if and only if $\chi_{E} \in
BV(\Omega)$ and $P(E;\Omega)$ coincides with $|D\chi_{E}|(\Omega)$, the total
variation in $\Omega$ of the distributional derivative of $\chi_{E}$.
Moreover,  a
generalised Gauss-Green formula holds:
\begin{equation}\nonumber
	{\int_{E}\mathrm{div}\varphi(x) \, dx
		=\int_{\Omega}\left\langle\nu_{E}(x),\varphi(x)\right\rangle \, d|D\chi_{E}|,
		\;\;\forall\,\varphi\in C_{c}^{1}(\Omega;\mathbb{R}^{N})},
\end{equation}
where $D\chi_{E}=\nu_{E}|D\chi_{E}|$ is the polar decomposition of $D\chi_{E}$.

\subsection{BH fields}\label{BH}

The space of fields with bounded hessian, introduced in \cite{D1, D4} and denoted by $BH(\Omega;\mathbb R^d)$, is defined as
\begin{align*}
BH(\Omega;\mathbb R^d):=&\left\{u \in W^{1,1}(\Omega;\mathbb R^d): D(\nabla u) \hbox{ is a bounded Radon measure }\right\} \\
=& \left\{u \in L^1(\Omega;\mathbb R^d): 
Du \in BV(\Omega;\mathbb R^{d\times N})\right\},
\end{align*}
the measure $D(\nabla u)$ is often written as $D^2u$.
The space $BH(\Omega;\mathbb R^d)$ is endowed with the norm
\begin{align*}
\|u\|_{BH(\Omega;\mathbb R^d)} :=\|u\|_{W^{1,1}(\Omega;\mathbb R^d)}
+ |D^2 u|(\Omega),
\end{align*}
where the latter term denotes the total variation of the Radon measure $D^2 u$.

Given $u \in BH(\Omega;\mathbb R^d)$, 
taking into account that $Du=\nabla u \in BV(\Omega;\mathbb R^d)$, the following decomposition holds, 
\begin{equation*}
D(\nabla u)= \nabla^2 u \mathcal L^N + D^s(\nabla u) =
\nabla^2 u \mathcal L^N + 
[\nabla u] \otimes \nu_{\nabla u} \mathcal H^{N-1}\lfloor S_{\nabla u} 
+ D^c(\nabla u).
\end{equation*}

Due to the symmetry of the distribution $D^2 u = D(\nabla u)$, in the previous expression,  $\nabla^2 u \mathcal L^N$ is a vector-valued measure, in the space of 
$d\times N \times N$  symmetric tensors $\mathbb R^{d \times N\times N}_s$, that is absolutely continuous with respect to $\mathcal L^N$ and with 
$$\displaystyle \nabla^2 u= \frac{\partial ^2 u_i}{\partial x_j \partial x_k}= \frac{\partial^2 u_i}{\partial x_k \partial x_j}, \; 1 \leq i \leq d, \; 1 \leq j,k \leq N,$$ 
$S_{\nabla u}$ denotes the singular set of $\nabla u$, $\nu_{\nabla u}$ is the normal to $S_{\nabla u}$ , $[\nabla u]:= (\nabla u)^+- (\nabla u)^-$
is the jump across $S_{\nabla u}$ and $D^c(\nabla u)$ is the Cantor part of $D(\nabla u)$, which is singular with respect to  $\mathcal L^N \lfloor \Omega+ \mathcal H^{N-1}\lfloor S_{\nabla u}$.
 
Again by the symmetry of $D^2u$, it follows that $[\nabla u] = \alpha  \otimes \nu_{\nabla u}$, for some function $\alpha$ which is integrable in $\Omega$ with respect to the measure $\mathcal H^{N-1}\lfloor S_{\nabla u}$, due to the fact that
$[\nabla u] \otimes \nu_{\nabla u} $ is a $d$-tuple of $N \times N$ symmetric rank-one matrices.

The space $SBH(\Omega;\mathbb R^d)$ consists of those functions $u \in BH(\Omega;\mathbb R^d)$ for which $D^c(\nabla u) = 0$,
that is,
$$SBH(\Omega;\mathbb R^d):= \left\{u \in L^1(\Omega;\mathbb R^d): 
Du \in SBV(\Omega;\mathbb R^{d\times N})\right\}.$$

In the next theorem we recall the approximate differentiability results that are valid in $BH$, see \cite[Theorem 2.1]{FHP}.

\begin{Theorem}\label{thm2.1FHP}
If $u \in BH(\Omega;\mathbb R^d )$, then
\begin{itemize}
\item[(i)] for $\mathcal L^N$ a.e. $x \in \Omega$
\begin{align*}
\lim_{\e \to 0^+}\frac{1}{\e^{N+2}}
\int_{Q(x,\e)}\left|u(y) - u(x)-\nabla u(x)(y - x) - 
\frac{1}{2}\nabla^ 2u(x)(y- x, y -x)\right| dy = 0,
\end{align*}
and
\begin{align*}
\lim_{\e\to 0^+}\frac{1}{\e^{N+1}}\int_{Q(x,\e)}
\left|\nabla u(y)-\nabla u(x)-\nabla ^2u(x)(y -x)\right| dy = 0;
\end{align*}
\item[(ii)] for $ {\mathcal H}^{N-1}$ a.e. $x \in S_{\nabla u}$ we have
\begin{align*}
\lim_{\e \to 0^+}\frac{1}{\e^{N+1}}\int_{Q^{\pm}_{\nu}(x,\e)}
\left|u(y) - u(x)-\nabla u^{\pm}(x)(y-x)\right|dy = 0,
\end{align*}
and
\begin{align*}
\lim_{\e \to 0^+}\frac{1}{\e^{N}}\int_{Q^{\pm}_{\nu}(x,\e)}
|\nabla u(y)-\nabla u^{\pm}(x)| dy = 0,
\end{align*}
where $Q^{\pm}_{\nu}(x,\e) = Q_{\nu}(x,\e) \cap 
\left\{y:\pm(y -x) \cdot \nu(x) > 0\right\}$. 
\end{itemize}
\end{Theorem}

Regarding embeddings of the space $BH(\Omega;\mathbb R^d)$, the following holds
(see \cite{D1}).
	
\begin{Proposition}\label{embedding}
	Let $\Omega \subset \mathbb R^N$ be an open, bounded set with Lipschitz boundary.
	Then
	$BH(\Omega;\mathbb R^d)\subset W^{1,p}(\Omega;\mathbb R^d)$
	with continuous embedding if $p \leq \frac{N}{N-1}$, the embedding is compact if $p < \frac{N}{N-1}$.
\end{Proposition}

The proof of the following interpolation inequality can be found in \cite[Theorem 2.2]{FLP2}.

\begin{Theorem}\label{interpolBH}[Interpolation inequality]
	Let $\Omega \subset \mathbb R^N$ be an open, bounded set with Lipschitz boundary. 
	Then, for every $\varepsilon > 0$, there is a constant $C = C(\varepsilon)$ such that
	$$\|\nabla u\|_{L^1(\Omega;\mathbb R^{d \times N})} \leq C\|u\|_{L^1(\Omega;\mathbb R^d)} + \varepsilon|D^2u|(\Omega),$$
	for all $u \in BH(\Omega;\mathbb R^d)$.
\end{Theorem}

\medskip

The next proposition states that, under the required hypotheses on $\Omega$, $BH$ functions can be extended to the whole space without charging the boundary, see \cite[Theorem 5.7]{H}.

\begin{Proposition}
	\label{thm5.7Hageryy}
	Let $\Omega \subset \mathbb R^N$ be an open, bounded set with Lipschitz boundary.  For any function $u$ in $BH(\Omega;\mathbb R^d)$ there exists an extension $E[u] \in BH(\mathbb R^N; \mathbb R^d)$ such that $|D(\nabla E[u])|(\partial \Omega)=0$.
\end{Proposition}

In some instances, as in Theorem \ref{cor5.8H} below, we will need a stronger notion of convergence for measures, known as area-strict convergence, written symbolically as (area) $\langle \cdot \rangle$-strict convergence. This notion is recalled in the following definition, see \cite{Rbook}.

\begin{Definition}\label{1100}
	Let $O \subset \mathbb R^N$ be a Borel set and, for \(n\in\mathbb N\),
	let $\mu_n=m_n^a\mathcal L^N+\mu_n^s\in\mathcal M(O;\mathbb
	R^{d})$  and  $\mu=m^a\mathcal L^N+\mu^s\in\mathcal M(O;\mathbb
	R^{d})$  be Radon measures on \(O\). We say that
	\begin{itemize}
\item[(i)] $\{\mu_n\}_{n\in\mathbb N}$ converges
		strictly to $\mu$ if $\mu_n\overset{*}{\rightharpoonup}\mu$
		and $|\mu_n|(O)\to|\mu|(O)$;
		\item[(ii)] $\{\mu_n\}_{n\in\mathbb N}$
		(area) $\langle\cdot\rangle$-strictly converges to $\mu$ if $\mu_n\overset{*}{\rightharpoonup}\mu$
		and $\langle\mu_n\rangle(O)\to\langle\mu\rangle(O)$,
		where for every measure $\lambda \in \mathcal M(O;\mathbb R^d)$,
		\begin{equation*}
			\langle \lambda\rangle(A):=\int_A \sqrt{1+|\lambda ^a(x)|^2}\,d x+|\lambda^s|(A),
			\quad\text{for every \(A\in \mathcal B(O)\)}.
		\end{equation*}
		\end{itemize}
\end{Definition}

\begin{Theorem}\label{cor5.8H}
Let $\Omega \subset \mathbb R^N$ be an open, bounded set with Lipschitz boundary. 
For any function $u\in BH(\Omega;\mathbb R^d)$, there exist smooth functions 
$\{u_n\} \subset C^{\infty}(\Omega;\mathbb R^d)$ such that $u_n \to u$ in $L^1(\Omega;\mathbb R^d)$, $\nabla u_n \to  \nabla u$ in
$L^1(\Omega;\mathbb R^{d \times N})$, $\nabla^2 u_n \overset{\ast}{\rightharpoonup} D(\nabla u)$ in $\mathcal M(\Omega;\Rb_s^{d \times N\times N})$ and $\nabla^2 u_n \mathcal L^N$ (area) $\langle \cdot \rangle$-strictly converges to $ D^2 u$, namely
$$
\int_\Omega \sqrt{1+ |\nabla^2 u_n(x)|^2} \, dx \to 
\int_\Omega \sqrt{1+ |\nabla^2 u(x)|^2} \, dx+ |D^s(\nabla u)|(\Omega).
$$
\end{Theorem}

The proof of this density result may be found in \cite[Corollary 5.8]{H}.

\medskip

For every $u \in BH(\Omega;\mathbb R^d)$, we define the spaces 
\begin{align*}
W^{2,1}_u(\Omega;\mathbb R^d):=\{v \in W^{2,1}(\Omega;\mathbb R^m): 
w \in BH(\mathbb R^N;\mathbb R^d) \hbox{ and } 
|D^2 w|(\partial \Omega)=0\}\end{align*}
and
\begin{align*}
BH_u(\Omega;\mathbb R^d):=\{v \in BH(\Omega;\mathbb R^d): 
w\in BH(\mathbb R^N;\mathbb R^d) \hbox{ and } |D^2w|(\partial \Omega)=0\},
\end{align*}
where $ w:=\begin{cases}
u-v, &\hbox{ in }\Omega\\
0, & \hbox{ in }\mathbb R^N \setminus \overline \Omega.
\end{cases}$

\medskip

Then, the following result is proved in \cite[Theorem 2.9]{FMZ}.

\begin{Lemma}\label{Lemma1KRBH}
	Let $\Omega\subset \mathbb R^N$ be a non-empty, bounded and open set such that $\mathcal L^N(\partial \Omega$)=0. For each $u \in BH(\Omega;\mathbb R^d)$ there exists $\{v_j\}\subset W^{2,1}_u(\Omega;\mathbb R^d)\cap C^\infty(\Omega;\mathbb R^d)$ such that $v_j \to u$ in $W^{1,1}(\Omega;\mathbb R^d)$, $\nabla^2 v_j \to D^2u$ (area) $\langle \cdot\rangle$-strictly in $\Omega$. If $u \in W^{2,1}(\Omega;\mathbb R^d)$,
	then we may, in addition, ensure that $v_j \to u$ strongly in $W^{2,1}(\Omega;\mathbb R^d)$.
\end{Lemma}

We end this subsection with a lemma, whose proof may be found in \cite{FM}, that will be used in Section \ref{mainsec}.

\begin{Lemma}\label{FM}
Let $\lambda$ be a non-negative Radon measure in $\mathbb R^N$. For $\lambda$ a.e. $x_0 \in \mathbb R^N$, for every $0 < \delta < 1$ and for every $\nu \in \mathbb \mathbb S^{N-1}$, the following holds
$$\limsup_{\varepsilon \to 0^+}\frac{\lambda(Q_\nu(x_0,\delta \varepsilon))}{\lambda(Q_\nu(x_0,\varepsilon))} \geq \delta^N,$$
so that
$$\lim_{\delta \to 1^-}\limsup_{\varepsilon \to 0^+}\frac{\lambda(Q_\nu(x_0,\delta \varepsilon))}{\lambda(Q_\nu(x_0,\varepsilon))} = 1.$$
\end{Lemma} 

\subsection{Notions of quasiconvexity}\label{defqcx}

 The concept of $2$-quasiconvexity, which we recall here, was introduced by Meyers in \cite{Mey}, extending Morrey's notion of quasiconvexity \cite{Mor} to integrands depending on second-order gradients.
\color{black}

\begin{Definition}\label{2qcx}
	A Borel measurable 
	function $f:\mathbb R^{d\times N\times N}_s\to \mathbb R$ is said to be 
	{\em 2-quasiconvex} if
	\begin{equation*}
		f(\xi)\leq\int_Q 
		f(\xi+\nabla^2 \varphi(x)) \, dx,
	\end{equation*}
	for every $\xi \in \mathbb R^{d\times N\times N}_s$ and for 
	every $\varphi \in C^{\infty}_{0}(Q;\mathbb R^d)$.  
\end{Definition}

\begin{Remark}\label{LDper}
	{\rm The above property is independent of the size, orientation and centre of the cube over which the integration is performed. Also, if $f$ is upper semi-continuous and locally bounded from above, using Fatou's Lemma and the density of smooth functions in $W^{2,1}(Q;\mathbb R^d)$, it follows that in Definition \ref{2qcx} $C^{\infty}_{0}(Q;\mathbb R^d)$ may be replaced by $W^{2,1}_{0}(Q;\mathbb R^d).$}
\end{Remark}

Given $f:\mathbb R^{d\times N\times N}_s\to \mathbb R$, the $2$-quasiconvex envelope of $f$, $Q^2f$, is defined by
\begin{equation}\label{2qcxenv}
	Q^2f(\xi) : = \inf \bigg\{\int_Q f(\xi+\nabla^2\varphi(x)) \, dx : \varphi \in C^{\infty}_{0}(Q;\mathbb R^d)\bigg\}.
\end{equation}
It is possible to show that,  if $f$ is upper semi-continuous,  $Q^2f$ is the greatest $2$-quasiconvex function that is less than or equal to $f$. Moreover, definition \eqref{2qcxenv} is independent of the domain, i.e.
\begin{align*}Q^2f(\xi) : = \inf \bigg\{\frac{1}{\L ^N(D)}\int_{D} f(\xi+\nabla^2\varphi(x)) \, dx : \varphi \in C^{\infty}_{0}(D;\mathbb R^d)\bigg\}
\end{align*}
whenever $D \subset \Rb^N$ is an open, bounded set with $\L ^N(\partial D) = 0$.

\medskip 

The proofs of the following results are based on the same arguments given in \cite[Remark 3.2 parts $(iv), (v)$]{CZ} and \cite[Propositions 2.6 and 2.7]{RZ}, see also \cite[Propositions 3.1 and 3.2]{BMZ} for a version in the framework of symmetrised quasiconvex envelopes.

\begin{Proposition}\label{SQfinfty=}
	Let $f:\{0,1\} \times \mathbb R^{d\times N \times N}_s\to [0, + \infty)$ be a continuous function as in \eqref{densityint} and satisfying \eqref{growthint} and \eqref{finfty}. Let $f^\infty$ and $Q^2f$ be its recession function and its $2$-quasiconvex envelope, with respect to the second variable, defined by \eqref{recS} and \eqref{2qcxenv}, respectively.
	Then
	\begin{equation*}
		Q^2(f^\infty)(q,\xi)= (Q^2f)^\infty(q,\xi) \;\;\;\hbox{ for every }(q,\xi) \in \{0,1\} \times \mathbb R^{d\times N \times N}_s.
	\end{equation*}
\end{Proposition}

\begin{Proposition}\label{propperH5}	
	Let $f:\{0,1\} \times \mathbb R^{d\times N \times N}_s\to [0, + \infty)$ be a continuous function as in \eqref{densityint}, satisfying \eqref{growthint} and \eqref{finfty}. 
	Then, there exist $\gamma\in [0,1)$ and $C>0$ such that
	$$ \displaystyle\left| (Q^2f)^\infty(q,\xi)- Q^2f(q, \xi)\right|
	\leq C\big( 1+|\xi|^{1-\gamma}\big), \;\;
	\forall\ (q,\xi)\in \{0,1\} \times \mathbb R^{d\times N \times N}_s.$$
\end{Proposition}

\section{Auxiliary results}\label{auxres}

The growth conditions \eqref{growthint}, as well as standard diagonalisation arguments, allow us to prove the following properties of the functional 
$\mathcal F (\chi,u;A)$ defined in \eqref{calFint}.

\begin{Proposition}\label{firstprop}
	Let $A \in \mathcal O(\Omega)$, $u \in BH(A;\mathbb R^d)$, $\chi \in BV(A;\{0,1\})$ and $\mathcal F(\chi,u;A)$ be given by \eqref{calFint}, where
	$F(\chi,u;A)$ is as in \eqref{Fint}. If $W_i$, $i = 0,1$, satisfy \eqref{growthint}, then
	\begin{itemize}
		\item[i)] there exists $C > 0$ such that 
		\begin{align*}C\left(|D^2u|(A) + |D\chi|(A)\right) \leq \mathcal F (\chi,u;A) \leq 
		C\left(\L^N(A) + |D^2u|(A) + |D\chi|(A)\right);
		\end{align*}
		\item[ii)] $\mathcal F (\chi,u;A)$ is always attained, that is, there exist sequences
		$\{u_n\} \subset  W^{2,1}(A;\mathbb{R}^{d})$ and
		$\{\chi_{n}\} \subset BV(A;\{0,1\})$ such that $u_{n}\to u$ in $W^{1,1}(A;\mathbb{R}^d)$,
		$\chi_{n}\to\chi$ in $L^1(A;\{0,1\})$ and
		$$\mathcal F(\chi,u;A) = \lim_{n\rightarrow+\infty}F(\chi_n,u_n;A);$$
		\item[iii)] if $\{u_n\} \subset BH(A;\mathbb{R}^{d})$ 
			and
		$\{\chi_{n}\} \subset BV(A;\{0,1\})$ are such that $u_{n}\to u$ in $W^{1,1}(A;\mathbb{R}^d)$ and 
		$\chi_{n}\to\chi$ in $L^1(A;\{0,1\})$, then
		$$\mathcal F(\chi,u;A) \leq \liminf_{n\to +\infty}\mathcal F(\chi_n, u_n;A);$$
		\item[iv)] $\mathcal F$ is local, that is, 
		\begin{align*}
		\mathcal F(\chi, u; A)= \mathcal F(\eta, v; A),
		\end{align*}
		for every $\chi, \eta \in BV(\Omega;\{0,1\})$ with  $\chi = \eta$, 
		$\mathcal L^N$ a.e. in  $A$, and for every $u,v \in BH(\Omega;\mathbb R^d)$ with $u=v$, $\mathcal L^N$ a.e. in  $A.$
	\end{itemize}
\end{Proposition}
\begin{proof}

i) The upper bound follows from the growth condition from above of $W_i$ \eqref{growthint}, 
$i  =0,1$ and by fixing, as test sequences for 
$\mathcal F (\chi,u;A)$,  $\chi_n = \chi$ and  $\{u_n\} \subset  W^{2,1}(A;\mathbb{R}^{d})$, such that
$u_{n}\to u$ in $W^{1,1}(A;\mathbb{R}^d)$ and $\displaystyle \lim_n|D^2u_n|(A) = |D^2u|(A)$. This choice is possible as the functional 
$u \mapsto |D^2u|(A)$ is lower semi-continuous with respect to the weak* convergence in $BH$, hence it coincides with its relaxation. Thus, the following
chain of inequalities yields the upper bound: 
\begin{align*}
\mathcal F (\chi,u;A) & \leq \liminf_{n \to +\infty}F(\chi,u_n;A) \\
& \leq \liminf_{n \to +\infty}\left[\int_A C\big(1 + |\nabla^2u_n(x)|\big) \, dx +
|D\chi|(A)\right]\\
& \leq \liminf_{n \to +\infty}C\left[\L^N(A) + |D^2u_n|(A) + |D\chi|(A)\right] \\
& \leq C\left[\L^N(A) + |D^2u|(A) + |D\chi|(A)\right].
\end{align*}
On the other hand, the lower bound is a consequence of the inequality from below in \eqref{growthint} and the lower semi-continuity of the total variation of Radon measures.

The conclusions in ii) and iii) follow by standard diagonalisation arguments.
As for the last statement, its proof is omitted since it is entirely similar to the one presented in \cite[Lemma 5.2]{FHP}.
\end{proof}

The same arguments which lead to the proof of \cite[Proposition 3.9]{BMZ}, replacing therein symmetrised gradients by hessians, enable us to obtain the following result.

\begin{Proposition}\label{newslicing}
Let $A \in \mathcal O (\Omega)$ and assume that $W_0, W_1$ satisfy the growth condition \eqref{growthint}. Let
$\{u_n\}$ and $\{v_n\}$ belong to $BH(A;\mathbb R^d)$ and $\{\chi_n\}, \{\theta_n\} \subset BV(A;\{0,1\})$ be sequences
	satisfying $u_{n} - v_n \to 0$ in $W^{1,1}(A;\mathbb{R}^{d})$,
	$\chi_{n} - \theta_n \to 0$ in $L^1(A;\{0,1\})$, 
	$\sup_n |(D^2)^su_n|(A) < + \infty$, $|D^2v_n| \overset{\ast}{\rightharpoonup} \mu$,
	$|D^2v_n|(A) \to \mu(A)$, $\sup_n |D\chi_n|(A) < + \infty$ and 
	$\sup_n |D\theta_n|(A) < + \infty$. Then there exist subsequences $\{v_{n_k}\}$ of $\{v_n\}$, 
	$\{\theta_{n_k}\}$ of $\{\theta_n\}$ and there exist sequences 
	$\{w_k\} \subset BH(A;\mathbb R^d)$, $\{\eta_k\} \subset BV(A;\{0,1\})$
	such that ${\rm supp}(w_k - v_{n_k}) \Subset A$, 
	${\rm supp}(\eta_k -\theta_{n_k}) \Subset  A$,
	$w_k - v_{n_k} \to 0$ in $W^{1,1}(A;\mathbb{R}^d)$,
	$\eta_{k} - \theta_{n_k} \to 0$ in $L^1(A;\{0,1\})$ and
	$$\limsup_{k\to +\infty}F(\eta_{k},w_k;A) \leq 
	\liminf_{n\to +\infty}F(\chi_n,u_n;A).$$
\end{Proposition}
\begin{proof}[Proof]
The proof follows as in \cite[Proposition 3.7]{BFT}, adjusted to take into account that we are now dealing with $BH$ fields, and also using the argument presented in the proof of  \cite[Proposition 3.9]{BMZ} 
which enables us to glue two characteristic functions of sets with finite perimeter in such a way that no new interface is created.
\end{proof}

\begin{Remark}\label{remreg}
{\rm Observe that the new sequence $\{w_k\}$ has the same regularity as the original sequences $\{u_n\}$ and $\{v_n\}$ as it is obtained through a convex combination of these ones using smooth cut-off functions.}
\end{Remark}

\begin{Proposition}\label{nestedsa}
	Assume that $W_0$ and $W_1$ are continuous functions satisfying \eqref{growthint}. 
	Let $u \in BH(\Omega;\mathbb R^d)$, $\chi \in BV(\Omega;\{0,1\})$ and $S, U, V \in \mathcal O(\Omega)$ be such that $S \subset \subset V \subset U.$ Then
	$$\mathcal{F}\left(\chi,u;U\right) \leq \mathcal{F}\left(\chi,u;V\right) 
	+ \mathcal{F}\left(\chi,u;U\setminus \overline S\right).$$
\end{Proposition}

\begin{proof}
	By Proposition \ref{firstprop}, ii), let 
	$\{v_n\} \subset W^{2,1}(V;\mathbb{R}^{d})$, 
	$\{w_n\} \subset W^{2,1}(U\setminus \overline S;\mathbb{R}^{d})$, 
	$\{\chi_n\} \subset BV(V;\{0,1\})$ and
	$\{\theta_n\} \subset BV(U\setminus \overline S;\{0,1\})$ be such that
	$v_n \to u$ in $W^{1,1}(V;\mathbb{R}^{d})$,
	$w_n \to u$ in $W^{1,1}(U\setminus \overline S;\mathbb{R}^{d})$,
	$\chi_{n}\to\chi$ in $L^1(V;\{0,1\})$
	$\theta_{n}\to\chi$ in 
	$L^1(U\setminus \overline S;\{0,1\})$ and
	\begin{align}\label{VUS1}
		\mathcal{F}\left(\chi,u;V\right) & = \lim_{n\rightarrow+\infty}F(\chi_n,v_n;V) \\
		\mathcal{F}\left(\chi,u;U\setminus \overline S\right) & = \lim_{n\rightarrow+\infty}F(\theta_n,w_n;U\setminus \overline S). \label{VUS2}
	\end{align}

Let $V_0 \in \mathcal O_{\infty}(\O)$ satisfy 
$S \subset \subset V_0 \subset \subset V$ and $|D^2u| (\partial V_0) = 0$,
$|D\chi|(\partial V_0) = 0$.

	Applying Proposition \ref{newslicing} to $\{v_n\}$ and $u$, $\{\chi_n\}$ and $\chi$, in $V_0$, we obtain sequences $\{\overline \chi_n\} \subset BV(V_0;\{0,1\})$  and $\{\overline v_n\} \subset W^{2,1}(V_0;\mathbb{R}^{d})$
	such that $\overline v_n = u$ near $\partial V_0$, $\overline v_n \to u$ in 
	$W^{1,1}(V_0;\mathbb{R}^{d})$, $\overline \chi_n = \chi$ near $\partial V_0$, $\overline \chi_n \to \chi$ in $L^{1}(V_0;\{0,1\})$  and
	\begin{equation}\label{sl1}
		\limsup_{n\to +\infty}F(\overline \chi_n,\overline v_n;V_0) \leq 
		\liminf_{n\to +\infty}F(\chi_n,v_n;V_0).
	\end{equation}
	A further application of Proposition \ref{newslicing}, this time to $\{w_n\}$ and $u$, $\{\theta_n\}$ and $\chi$, in $U \setminus \overline V_0$, yields sequences $\{\overline \theta_n\} \subset BV(U \setminus \overline V_0;\{0,1\})$,  
	$\{\overline w_n\} \subset W^{2,1}(U \setminus \overline V_0;\mathbb{R}^d)$,
	such that $\overline w_n = u$ near $\partial V_0$, $\overline w_n \to u$ in 
	$W^{1,1}(U \setminus \overline V_0;\mathbb{R}^d)$, $\overline \theta_n = \chi$ near $\partial V_0$, $\overline \theta_n \to \chi$ in 
		$L^{1}(U \setminus \overline V_0;\{0,1\})$ and
	\begin{equation}\label{sl2}
		\limsup_{n\to +\infty}F(\overline \theta_n,\overline w_n;U \setminus \overline V_0) \leq \liminf_{n\to +\infty}F(\theta_n,w_n;U \setminus \overline V_0).
	\end{equation}
	Define 
	$$z_n : = \begin{cases}
		\overline v_n,& {\rm in } \; V_0\\
		\overline w_n,& {\rm in } \; U \setminus V_0,
	\end{cases}$$
	notice that, by the properties of $\{\overline v_n\}$ and $\{\overline w_n\}$,
	$\{z_n\} \subset W^{2,1}(U;\mathbb{R}^{d})$ and 
	$z_n \to u$ in $W^{1,1}(U;\mathbb{R}^{d})$.

	We must now build a transition sequence $\{\eta_n\}$ between $\{\overline\chi_n\}$ and $\{\overline \theta_n\}$, in such a way that an upper bound for the total variation of $\eta_n$ is obtained.
	In order to connect these functions without adding more interfaces, we apply the argument used in the proof of \cite[Proposition 3.9]{BMZ}. 
	
In this way we show the existence of a set with locally Lipschitz boundary
$V_{\rho_0} \subset V$ so that $\overline w_n = u$ in $V_{\rho_0} \setminus \overline V_0$,
\begin{equation}\label{rho1}
\int_{V_{\rho_0} \setminus \overline V_0}\big(1 + |\nabla^2u(x)|\big) \, dx = O(\rho_0),
\; \; |D\chi|(\partial V_{\rho_0}) =0,
\end{equation}
and
\begin{equation}\label{rho2}
\lim_{n \to + \infty}\int_{\partial V_{\rho_0}}|\overline \chi_n(x) - \theta_n(x)| \, d \mathcal H^{N-1}(x) = 0.
\end{equation}
Defining
\begin{equation*}
\eta_{n} :=\begin{cases}
\overline \chi_n, &\hbox{ in } V_{\rho_0}\\
\overline \theta_n, &\hbox{ in } U \setminus V_{\rho_0},
\end{cases}
\end{equation*}
it turns out that $\{\eta_{n}\} \subset BV(U;\{0,1\})$ and
$\eta_{n}\to\chi$ in $L^1(U;\{0,1\})$
so that $\{\eta_{n}\}$
and $\{z_{n}\}$ are admissible for $\mathcal{F}\left(\chi,u;U\right)$.
Hence, by \eqref{growthint}, \eqref{rho1}, \eqref{rho2}, \eqref{sl1}, \eqref{sl2}, \eqref{VUS1}, \eqref{VUS2}
we obtain
\begin{align*}
\mathcal{F}\left(\chi,u;U\right) &\leq 
\liminf_{n\to +\infty}F(\eta_{n},z_{n};U) \\
& \leq \limsup_{n\to +\infty}F(\overline\chi_n,\overline v_n; V_0) + 
\limsup_{n\to +\infty}F(\overline\theta_n,\overline w_n; U \setminus \overline V_{0})
+ \int_{V_{\rho_0}\setminus \overline V_0}C \big(1 + |\nabla^2 u(x)|\big) \, dx 
\end{align*}
\begin{align*}
& \hspace{3cm} 
+ \limsup_{n\to +\infty}\int_{\partial V_{\rho_0} }|\chi_n(x) - \theta_n(x)| 
\, d  {\mathcal H}^{N-1}(x) 
+ \limsup_{n\to +\infty} |D\overline \chi_n|(V_{\rho_0}\setminus V_0) \\
& \leq \liminf_{n\to +\infty}F(\chi_n,v_n; V_0) +
\liminf_{n\to +\infty}F(\theta_n,w_n; U \setminus \overline V_{0}) + O(\rho_0) +
|D \chi|(V_{\rho_0}\setminus V_0) \\
& \leq \mathcal{F}\left(\chi,u;V\right) + 
\mathcal{F}\left(\chi,u;U\setminus \overline S\right) + O(\rho_0) +
|D \chi|(V_{\rho_0}\setminus V_0), 
\end{align*}
so the result follows by letting $\rho_0 \to 0^+$, taking into account 
that $|D\chi|(\partial V_0)=0$ and \eqref{rho1}.
\end{proof}

The proof of the following result is identical to the one of \cite[Proposition 3.8]{BMZ}, hence it is omitted. 

\begin{Proposition}\label{traceRm}
	Let $W_0$ and $W_1$ be continuous functions satisfying \eqref{growthint}. 
	For every $u \in BH(\Omega;\mathbb R^d)$, $\chi \in BV(\Omega;\{0,1\})$, $\mathcal{F}\left(\chi,u;\cdot\right)$ is the restriction to $\mathcal O(\O)$ of a Radon measure.
\end{Proposition}

\begin{Remark}\label{A1A2}
{\rm Due to Propositions \ref{traceRm} and \ref{firstprop} i), given any $(\chi,u) \in BV(\Omega;\{0,1\}) \times BH(\Omega;\mathbb R^d)$ and any open sets $A_1 \subset \subset A_2 \subseteq \Omega$, it follows that 
\begin{align*}
\mathcal F(\chi,u;A_2) \leq  \, \mathcal F(\chi,u;A_1) 
+C\Big(\mathcal L^N(A_2 \setminus A_1) + |D^2 u|(A_2 \setminus A_1)
+ |D\chi|(A_2 \setminus A_1)\Big).
\end{align*}
Indeed, for $\varepsilon > 0$ small enough, let $A_{\varepsilon} : = \left\{x \in A_1 : {\rm dist}(x, \partial A_1) > \varepsilon\right\}$ and notice that $A_2$ is covered by the union of the two open sets $A_1$ and $A_2 \setminus \overline{A_{\varepsilon}}.$ Thus, by Propositions \ref{traceRm} and \ref{firstprop} i), we have
\begin{align*}
\mathcal F(\chi,u;A_2) \leq & \, \mathcal F(\chi,u;A_1) + 
\mathcal F(\chi,u;A_2 \setminus \overline{A_{\varepsilon}}) \\
\leq & \, \mathcal F(\chi,u;A_1) +
C\Big(\mathcal L^N(A_2 \setminus  \overline{A_{\varepsilon}}) 
+ |D^2 u|(A_2 \setminus\overline{A_{\varepsilon}}) +
|D\chi|(A_2 \setminus\overline{A_{\varepsilon}})\Big).
\end{align*}
To conclude the result it suffices to let $\varepsilon \to 0^+$.
}
\end{Remark}

\smallskip

The following proposition, whose proof is standard (cf. for instance \cite[Lemma 3.1]{RZ} or \cite[Proposition 2.14]{CZasy}), relying on the relaxation result \cite[Theorems 1.1 and 3.6, together with Lemma 3.1]{BFL} allows us to assume without loss of generality that $f$ is $2$-quasiconvex with respect to the second variable. 

\begin{Proposition}\label{Frel=Frel**}
	Let $W_0$ and $W_1$ be continuous functions satisfying \eqref{growthint} and consider the functional 
	\linebreak \hfill
	$F:BV(\Omega;\{0,1\})\times BH(\Omega;\mathbb R^d)\times \mathcal O(\Omega) \to
	[0,+\infty]$ 
	defined in \eqref{Fint}.
	Consider furthermore the relaxed functionals given in \eqref{calFint} and
	\begin{align*}
		\mathcal F_{Q^2f}(\chi,u;A) & :=\inf\Big\{\liminf_{n\to +\infty}
		\left[\int_A  Q^2f(\chi_n(x),\nabla^2 u_n(x)) \, dx + |D \chi_n|(A)\right] : 
		 \{\chi_n\} \subset BV(A;\{0,1\}), \nonumber \\ 
		&\hspace{2cm}\{u_n\} \subset W^{2,1}(A;\mathbb R^d),  \chi_{n} \to \chi \text{ in } L^1(A;\{0,1\}),
		u_{n}\to u\text{ in } W^{1,1}(A;\mathbb{R}^d)\Big\}.
	\end{align*}
	
	Then, ${\mathcal F}(\cdot,\cdot;\cdot)$ coincides with 
	${\mathcal F}_{Q^2f}(\cdot, \cdot; \cdot)$ in 
	$BV(\Omega;\{0,1\})\times BH(\Omega;\mathbb R^d)\times \mathcal O(\Omega)$.
\end{Proposition}
\begin{proof}
	Relying on the results in \cite{BFL}, we begin by showing that, for every $A\in \mathcal O(\Omega)$ and every
	$(\chi,u) \in BV(A;\{0,1\}) \times W^{2,1}(A;\mathbb R^d)$,
	\begin{align}\label{repweak}
		&\inf\left\{\liminf_{n\to +\infty}\int_{A}f(\chi_n(x), \nabla^2 u_n(x)) \, dx
		: u_n \rightharpoonup u \hbox{ in } W^{2,1}(A;\mathbb R^d) , \chi_n \to \chi \in L^1(A;\{0,1\})\right\} \nonumber \\
		& \hspace{2cm}	=\int_{A} Q^2f(\chi(x), \nabla^2 u(x)) \,dx. 
	\end{align}
	Indeed, we apply the results in \cite{BFL} to the fields $(U,V)$, with $U=\chi \in L^1(A;\{0,1\})$ and $V=\nabla^2 u  \in  L^1(A; \mathbb R_{s}^{d \times N\times N})$, where the operator $\mathcal A$ appearing in \cite{BFL} acts on $V$ as follows
	$$
	{\mathcal A} V :=\left(
	\frac{\partial V_{jk}}{\partial x_i}
	- \frac{\partial V_{ji}}{\partial x_k}\right)_{1\leq i,j,k\leq N}.
	$$
	
	It was shown in \cite[Example 3.10(d)]{FMAq}, that 
	the fields $V$ which belong to the kernel of $\mathcal A$ are of the form
	$V = \nabla^2 u$ for some $u \in W^{2,1}(A;\mathbb R^{d})$.
Consequently, the $\mathcal A$-quasiconvex envelope of $f$ with respect to the second variable $\nabla^2u$ coincides with the 2-quasiconvexification of $f$ with respect to the same variable, i.e.
	\begin{align}\label{QA=Q2}
		Q_{\mathcal A} f(\chi,\nabla^2 u) = Q^2 f(\chi, \nabla^2 u),
		\forall (\chi,u) \in BV(A;\{0,1\}) \times BH(A;\mathbb R^d).
	\end{align}
	For a more general proof of this equality we refer to \cite{Ra}.
	
	Hence the lower bound for the infimum given in \eqref{repweak} follows directly by \cite[Theorems 1.1 and 3.6]{BFL}, by means of equality \eqref{QA=Q2}.
	The reverse inequality can be deduced by \cite[Theorem 1.3]{BFL} applied to the Carath\'eodory function $h(x, \xi):= f(\chi(x), \xi)$, and taking into account that
	\begin{align}\label{ineqweak}
		&\inf\left\{\liminf_{n\to +\infty}\int_{A}f(\chi_n(x), \nabla^2 u_n(x)) \, dx
		: u_n \rightharpoonup u \hbox{ in } W^{2,1}(A;\mathbb R^d) , \chi_n \to \chi \in L^1(A;\{0,1\})\right\}\nonumber\\
		& \leq \inf\left\{\liminf_{n\to +\infty}\int_{A}f(\chi(x), \nabla^2 u_n(x)) \, dx
		: u_n  \rightharpoonup u \hbox{ in } W^{2,1}(A;\mathbb R^d) \right\}.
	\end{align}
\end{proof}

\begin{Remark}\label{remrecession}{\rm
	In the sequel we rely on the result of Proposition \ref{Frel=Frel**} and assume
	that $f$ is  $2$-quasiconvex with respect to the second variable. Together with \eqref{growthint} and the continuity of $f$,  this entails the Lipschitz continuity of $f$ with respect to the second variable (see  \cite[Lemma 2.8 and Lemma 2.11]{H})  which, in turn, leads to the definition of the (strong) recession function considered in \eqref{recS}, see \cite{Rbook}.}
	\end{Remark}

We end this section with an analogue of Theorem \ref{main} in the case of $p$-growth ($p > 1$) conditions for the density functions.
	
\begin{Proposition}\label{casep}
Let $f:\{0,1\} \times \mathbb R^{d\times N \times N}_s\to [0, + \infty)$ be a continuous function as in \eqref{densityint}, where $W_0$ and $W_1$ satisfy \eqref{growthintp}, and consider 
$F:BV(\Omega;\{0,1\})\times BH(\Omega;\mathbb R^d)\times \mathcal O(\Omega) \to [0,+\infty]$ defined in \eqref{Fint}. Then
\begin{align}\label{repweakp}
\mathcal{F}\left(\chi,u;\Omega\right)   &  
=\inf\bigg\{ \liminf_{n\rightarrow +\infty} F(\chi_n,u_n;\Omega): 
\{\chi_n\} \subset BV(\Omega;\{0,1\}), \{u_n\} \subset W^{2,p}(\Omega;\mathbb{R}^{d}), \nonumber\\
& \hspace{4cm} \chi_{n}\to\chi\text{ in } L^1(\Omega;\{0,1\}),
u_{n}\to u\text{ in } W^{1,1}(\Omega;\mathbb{R}^{d})
  \bigg\}\nonumber\\
&= \begin{cases}
\displaystyle \int_\Omega Q^2f(\chi(x),\nabla^2 u(x)) \,dx + |D \chi|(\Omega), &\hbox{ if } 
(\chi, u)\in BV(\Omega;\{0,1\})\times W^{2,p}(\Omega;\mathbb R^d), \\
+\infty, & \hbox{ otherwise},
\end{cases}
\end{align}
where $\mathcal F(\chi,u;\Omega)$ is given in \eqref{calFint}.
\end{Proposition}
\begin{proof}
The proof of \eqref{repweakp} is very similar to that of \eqref{repweak}. Indeed, the arguments that were used to obtain the previous result also hold in the case  $p>1$. On the other hand, the coercivity condition in \eqref{growthintp} and classical embedding theorems allow us to consider indifferently weak convergence in $W^{2,p}$ or strong convergence in $W^{1,1}$. More precisely, we observe that the fact that the expression in \eqref{repweakp} is a lower bound for $\mathcal{F}\left(\chi,u;\Omega\right)$ is also a consequence of  the lower semi-continuity of $|D \chi|(\Omega)$ with respect to the $BV$-weak* convergence, and the opposite inequality can be proved as \eqref{ineqweak}, once $\chi$ is fixed.
\end{proof}

\section{Proof of the main result}\label{mainsec}
\subsection{The bulk and Cantor Terms}\label{Cantor}

By Proposition \ref{traceRm}, for every 
$(\chi, u)\in BV(\Omega;\{0,1\})\times BH(\Omega;\mathbb R^d)$, 
$\mathcal F(\chi,u,;\cdot)$ in \eqref{calFint} is the restriction to $\mathcal O(\Omega)$ of a Radon measure that we denote by $\mu$. 
	
By i) of Proposition \ref{firstprop}, $\mu$ can be decomposed as 
$$
\mu=\mu^a\mathcal L^N+ \mu^c+ \mu^j,
$$
where $\mu^j \ll |D^j(\nabla u)|+ |D \chi|$ and $\mu^c \ll |D^c(\nabla u)|$.

This subsection is devoted to the characterisation of the density $\mu^a$ and the measure $\mu^c$, that is, to the identification of the density of $\mu$ with respect to the measure $\lambda_u$ given by
\begin{equation}\label{lambdau}
\lambda_u:= \mathcal L^N + |D^c (\nabla u)|.
\end{equation} 

To this end, we start by observing that, by virtue of Proposition \ref{Frel=Frel**}, there is no loss of generality in assuming that $f$ is $2$-quasiconvex with respect to the second variable. 

\begin{Proposition}\label{lbCantor}
	Let $u \in BH(\Omega;\mathbb R^d)$, $\chi \in BV(\Omega;\{0,1\})$ and let $W_0$ and $W_1$ be continuous functions satisfying \eqref{growthint}. 
	Assume that $f$, given by \eqref{densityint}, is $2$-quasiconvex with respect to the second variable.
	Then,  for $\mathcal L^N$ a.e. $x_0\in \Omega$,
	$$
	\dfrac{d\mathcal F(\chi, u; \cdot)}{d \mathcal L^N}(x_0)\geq f(\chi(x_0), \nabla^2 u(x_0)),
	$$
	and, for $|D^c(\nabla u)|$ a.e. $x_0 \in \O$, 
	$$  \dfrac{d\mathcal F(\chi, u; \cdot)}{d |D^c(\nabla u)|}(x_0) \geq 
	f^{\infty}\left(\chi(x_0), \frac{d D^c(\nabla u)}{d |D^c(\nabla u)|}(x_0)\right).$$
\end{Proposition}
\begin{proof}

Let  $x_0 \in \O$ be a point satisfying
\begin{equation*}
\frac{d |(D^j(\nabla u)|}{d\lambda_u}(x_0)= 0, \; \; 
\frac{d |D\chi|}{d\lambda_u}(x_0)=0
\end{equation*}
and
\begin{equation}\label{RNdev}
\dfrac{d\mathcal F(\chi, u; \cdot)}{d\lambda_u}(x_0) = \frac{d \mu}{d \lambda_u}(x_0)
= \lim_{\e \to 0^+}\frac{\mu(Q(x_0,\e))}{\lambda_u(Q(x_0,\e))} 
\; \; \mbox{exists and is finite},
\end{equation}
the above properties hold for $\lambda_u$ a.e. point in $\Omega$.
In particular, \eqref{RNdev} tells us that, for $\mathcal L^N$ a.e. $x_0$, we have
\begin{equation}\label{RNdevleb}
\frac{d \mu}{d \lambda_u}(x_0)
= \lim_{\e \to 0^+}\frac{\mu(Q(x_0,\e))}{\e^N} =  \frac{d \mu}{d \mathcal L^N}(x_0)
\; \; \mbox{exists and is finite},
\end{equation}
and, for $|D^c(\nabla u)|$ a.e. $x_0$,
\begin{equation}\label{RNdevCantor}
\frac{d \mu}{d \lambda_u}(x_0) =
\lim_{\e \to 0^+}\frac{\mu(Q(x_0,\e))}{|D^c(\nabla u)|(Q(x_0,\e))} = \frac{d \mu}{d |D^c(\nabla u)|}(x_0)
\; \; \mbox{exists and is finite}.
\end{equation}
	
Furthermore, we choose $x_0$ to be a point of approximate continuity for $u$, for $\nabla u$, for $\nabla^2 u$ and for $\chi$, in particular
\begin{equation}\label{chiad}
\lim_{\e\rightarrow 0^+}\frac{1}{\e^{N}}
\int_{Q\left(x_0,\e\right)}\left\vert \chi(x)-\chi(x_0)\right\vert \,dx = 0.
\end{equation}
In view of \cite[(3.79)]{AFP}, the above properties hold for $\lambda_u$ a.e. $x_0 \in \Omega$.
	
In addition, recall that $\mathcal L^N$ a.e. $x_0$ also satisfies
\begin{equation*}
\lim_{\varepsilon \to 0^+} 
\frac{|D^s(\nabla u)|( Q(x_0,\varepsilon))}{\varepsilon^N}=0,
\end{equation*}
and that
\begin{equation}\label{derRNlambdau}
\frac{d |D(\nabla u)|}{d\lambda_u} (x_0)
= \lim_{\varepsilon \to 0^+}
\frac{|D (\nabla u)|(Q(x_0,\varepsilon))}{\lambda_u(Q(x_0,\varepsilon))}<+\infty
\end{equation}
holds for $\lambda_u$ a.e. $x_0\in \Omega$.

We define
$$f_0(\xi) = f(0, \xi) \mbox{ and } f_1(\xi) = f(1, \xi), \forall \xi \in \Rb^{d\times N \times N}_s$$
and, for $i=0,1$, we consider the auxiliary functionals
\begin{align}\label{auxf}
\mathcal F_i(u; A) &:= \inf \Big \{\liminf_{n \to + \infty} 
\int_A f_i (\nabla ^2 u_n(x)) \, dx :  \{u_n\} \subset W^{2, 1}(A; \Rb^d), u_n \to u \; {\rm in} \; W^{1,1}(A;\Rb^d)\Big \}. 
\end{align}
	From Theorem 3.1 in \cite{H}, we infer that 
	$\mathcal F_i(u;\cdot)$, $i=0,1$,  
	are the restriction to $\mathcal O(\Omega)$ of Radon measures whose densities with respect to $\mathcal L^N$ and $ |D^c(\nabla u)|$  are given, respectively, by 
	\begin{equation}\label{CdensFibulk}
	\frac{d \mathcal F_i(u;\cdot)}{d\mathcal L^N }(x_0) 
	= f_{i} \big(\nabla^2 u(x_0)\big) =
	f\big(i,\nabla^2 u(x_0)\big)
\end{equation}
for $\mathcal L^N$ a.e. $x_0 \in \Omega$, 
and
	\begin{equation}\label{CdensFi}
		\frac{d \mathcal F_i(u;\cdot)}{d|D^c (\nabla u)|}(x_0) 
		= f_{i}^\infty \Big(\frac{d D^c(\nabla u) }{d|D^c (\nabla u)|}(x_0)\Big) =
		f^\infty \Big(i,\frac{d D^c(\nabla  u)}{d|D^c (\nabla u)|}(x_0)\Big)
	\end{equation}
	for $|D^c(\nabla u) |$ a.e. $x_0 \in \Omega$.
	The point $x_0$ is selected so that it also satisfies these properties.

	In what follows we assume, without loss of generality, that $\chi(x_0)=1$, the case $\chi(x_0)=0$ can be treated similarly. Bearing this choice in mind we work with the 
	functional \eqref{auxf}, and we will make use of \eqref{CdensFibulk}, \eqref{CdensFi}, in the case $i=1$.
	
	Selecting the sequence $\e_k \to 0^+$ in such a way that 
	$\mu(\partial Q(x_0,\e_k)) = 0$ and $Q(x_0,\e_k) \subset \Omega$, we have
	\begin{align*}
		& \lim_{k \to +\infty}\frac{\mu (Q(x_0,\e_k))}{\lambda_u(Q(x_0,\e_k))} \\
		&=  \lim_{k,n}\frac{1}{\lambda_u (Q(x_0,\e_k))}\left[\int_{Q(x_0,\e_k)}
		f(\chi_n(x),\nabla^2 u_n(x))\,dx  +  |D \chi_n|(Q(x_0,\e_k))\right]
	\end{align*}
	where, due to ii) of Proposition \ref{firstprop}, 
	$\{\chi_n\} \subset BV(Q(x_0,\e_k);\{0,1\})$ and  $\{u_n\} \subset W^{2,1}(Q(x_0,\e_k);\Rb^d)$, sa\-tis\-fy $\chi_n \to \chi$ in $L^1(Q(x_0,\e_k);\{0,1\})$, $u_n \to u$ in $W^{1,1}(Q(x_0,\e_k);\Rb^d)$, and are recovery sequences for $\mathcal F(\chi, u; Q(x_0,\varepsilon_k))$.
	
	Taking into account that we are searching for a lower bound for the density of $\mu$ with respect to $\lambda_u$ at $ x_0$,
	we neglect the perimeter term $|D \chi_n|(Q(x_0,\e_k))$ and obtain
	\begin{align}\label{lbmuc}
		+\infty & > \lim_{k \to +\infty}\frac{\mu (Q(x_0,\e_k))}{\lambda_u(Q(x_0,\e_k))} \nonumber\\
		& \geq
		\liminf_{k,n}\frac{1}{\lambda_u(Q(x_0,\e_k))}
		\int_{Q(x_0,\e_k)}f(\chi_n(x),\nabla^2 u_n(x))\,dx  
		\nonumber\\
		& \geq
		\liminf_{k,n}\frac{1}{\lambda_u(Q(x_0,\e_k))}
		\int_{Q(x_0,\e_k)}f_1(\nabla^2 u_n(x))\,dx \nonumber\\
		& \hspace{3cm} + \limsup_{k,n}\frac{1}{\lambda_u(Q(x_0,\e_k))}
		\int_{Q(x_0,\e_k)}\big[f(\chi_n(x),\nabla^2 u_n(x)) - f(1,\nabla^2 u_n(x))\big]\,dx \nonumber\\
		& \geq \liminf_k \frac{\mathcal F_1(u;Q(x_0,\e_k))}{\lambda_u(Q(x_0,\e_k))}+ \limsup_{k,n}I_{k,n} \nonumber\\
		& \geq \frac{d \mathcal F_{1}(u;\cdot)}{d\lambda_u}(x_0) + \limsup_{k,n}I_{k,n} 
	\end{align}
	where
	$$I_{k,n} = \frac{1}{\lambda_u(Q(x_0,\e_k))}
	\int_{Q(x_0,\e_k)}\big[f(\chi_n(x),\nabla^2 u_n(x)) - f(1,\nabla^2 u_n(x))\big]\,dx.$$
	In view of \eqref{auxf}, \eqref{CdensFibulk} and \eqref{CdensFi},  it suffices to show that 
	$\displaystyle \limsup_{k,n}|I_{k,n}|=0$. Indeed, if this claim holds, then,  
	for $\mathcal L^N$ a.e. $x_0 \in \Omega$, we conclude that
	$$\frac{d \mu}{d \lambda_u}(x_0)= \frac{d \mu}{d \mathcal L^N}(x_0)\geq f_1(\nabla^2 u(x_0)) =
	f\big(\chi(x_0), \nabla^2u(x_0)\big)$$ 
	and, for $|D^c(\nabla u) |$ a.e. $x_0 \in \Omega$, we deduce that
	$$\frac{d \mu}{d \lambda_u}(x_0)=\frac{d \mu}{d |D^c(\nabla u)|}(x_0)\geq  
	f_1^\infty\left(\frac{dD^c (\nabla u)}{d |D^c (\nabla u)|}(x_0)\right) = 
	f^\infty\left(\chi(x_0),\frac{dD^c (\nabla u)}{d |D^c (\nabla u)|}(x_0)\right).$$

	It remains to prove the claim. Changing variables we obtain
	\begin{align}\label{Ikn}
		\left|I_{k,n}\right| &= \left|\frac{\e_k^N}{\lambda_u(Q(x_0,\e_k))}
		\int_Q \Big(f(\chi_n(x_0+\e_k y),\nabla^2 u_n(x_0+\e_k y)) - 
		f(1, \nabla^ 2 u_n(x_0+\e_k y))\Big)\,dy\right| \nonumber \\
		&= \left|\delta_k
		\int_Q (f(\chi_{n,k}(y) + 1,\nabla^2 u_{n,k}(y)) - 
		f(1, \nabla^2 u_{n,k}(y)))\,dy\right|,
	\end{align}
	where 
	$$\displaystyle \delta_k : = \frac{\e_k^N}{\lambda_u(Q(x_0,\e_k))}, \; \;
	\chi_{n,k}(y) := \chi_n(x_0+\e_k y) - 1, \; \; 
	u_{n,k}(y) := \frac{u_n(x_0+\e_k y)}{\e_k^2 }.$$
	By \eqref{chiad} it follows that 
	\begin{equation}\label{chiconv}
	\displaystyle \lim_{k,n}\|\chi_{n,k}\|_{L^1(Q)} = 0.
	\end{equation}
	On the other hand,
	due to \eqref{RNdevleb}, $\displaystyle \lim_k \delta_k = 1$, for $\mathcal L^N$ a.e. $x_0 \in \Omega$,
	whereas, by \eqref{RNdevCantor}, $\displaystyle \lim_k \delta_k = 0$, for $|D^c(\nabla u)|$ a.e. $x_0$.

	We argue as in \cite[Proposition 4.4 or 4.1]{BMZ}. Using \eqref{G}, \eqref{chiconv}, the growth condition from below on $f$ and \eqref{lbmuc}, we have from \eqref{Ikn},
	\begin{align}\label{Ikn2}
		\limsup_{k,n}\left|I_{k,n}\right| \leq 
		\limsup_{k,n} C \delta_k 
		\int_Q |\chi_{n,k}(y)| \, \big(1+|\nabla^2 u_{n,k}(y)|\big) \, dy <+\infty.
	\end{align}

Using a diagonalisation argument, let $\chi_k := \chi_{n(k),k}$, $u_k := u_{n(k),k}$ be such that
$\chi_k \to 0$ in $L^1(Q)$ and 
$$\limsup_{k,n} C \delta_k 
	\int_Q |\chi_{n,k}(y)| \, \big(1+|\nabla^2 u_{n,k}(y)|\big) \, dy 
	= \lim_k C\delta_k \int_Q |\chi_k(y)| \big(1+ |\nabla^2 u_k(y)|\big) \,dy.$$
	
Therefore, the sequence $\{\delta_k |\chi_k \,(1+ \nabla^2u_{k})|\}$ is bounded in $L^1(Q;\Rb^{d\times N \times N}_s)$
so, by Chacon's Biting Lemma, there exists a subsequence (not relabelled) and there exist
	sets $O_r \subset Q$ such that 
	$\displaystyle \lim_{r \to + \infty}{\mathcal L}^N(O_r) = 0$ and the sequence 
	$\{\delta_k |\chi_k \, (1+ \nabla^2 u_k)|\}$ is equiintegrable in $Q\setminus O_r$, 
	for any $r \in \mathbb N$.
	
This allows us to ensure that 
for any $j \in \mathbb N$ there exist $k(j), r(j) \in \mathbb N$ such that
\begin{equation}\label{Ikn4}
\delta_k(j)\int_{O_{r(j)}} |\chi_{k(j)}(y)| \, |1+\nabla^2 u_{k(j)}(y)| \, dy \leq \frac{1}{j}
	\end{equation}
and that, for any $\e > 0$, 
\begin{equation}\label{Ikn5}
\delta_{k(j)}\int_{Q \setminus O_{r(j)}} |\chi_{k(j)}(y)| \, |1 +\nabla^2 u_{k(j)}(y)| \, dy < \e,
\end{equation}
provided $j$ is large enough.

The conclusion now follows from \eqref{Ikn2}, \eqref{Ikn4} and \eqref{Ikn5}.
\end{proof}

\begin{Proposition}\label{ub}
	Let $u \in BH(\Omega;\mathbb R^d)$, $\chi \in BV(\Omega;\{0,1\})$ and let $W_0$ and $W_1$ be continuous functions satisfying \eqref{growthint}. 
	Assume that $f$ given by \eqref{densityint} is $2$-quasiconvex in the second variable and
	let $\lambda_u$ be as in \eqref{lambdau}. Then, for $\lambda_u$ a.e. $x_0 \in \O$,
$\dfrac{d\mathcal F(\chi, u; \cdot)}{d \lambda_u} (x_0)$ exists and is finite and we have
\begin{equation*}
\frac{d \mu}{d \mathcal L^N}(x_0)=\frac{d \mathcal F (\chi, u;\cdot)}{d \mathcal L^N}(x_0)\leq f(\chi(x_0),\nabla^2 u(x_0)), \; \mbox{ for }
	\mathcal L^N \mbox{ a.e } x_0 \in \Omega
	\end{equation*}
	and
\begin{equation*}
\frac{d \mu}{d |D^c(\nabla u)|}(x_0) = \frac{d \mathcal F (\chi, u;\cdot)}{d |D^c(\nabla u)|}(x_0)\leq f^\infty\left(\chi(x_0),\frac{d D^c(\nabla u)}{d |D^c(\nabla u)|}(x_0)\right), \; \mbox{ for }
	|D^c(\nabla u)| \mbox{ a.e } x_0 \in \Omega.
\end{equation*}
\end{Proposition}
\begin{proof}
	Let $x_0 \in \Omega$ be a point such that
	$\displaystyle	\frac{d \mu}{d \lambda_u}(x_0) = \dfrac{d\mathcal F(\chi, u; \cdot)}{d \lambda_u} (x_0)$ exists and is finite.
	We also select $x_0$ to be a point of approximate continuity for $u$, for $\nabla u$, for $\nabla^2 u$ and for $\chi$, in particular so that
	\eqref{chiad} and \eqref{derRNlambdau} hold and, in addition,
	\begin{equation}\label{Dchi}
		\lim_{\e \to 0^+}\frac{|D\chi|(Q(x_0,\e))}{\lambda_u(Q(x_0,\e))} = 0.
	\end{equation}
	As observed  in the proof of Proposition \ref{lbCantor} these properties are satisfied for $\lambda_u$ a.e. $x_0 \in \Omega$.
	
	Assuming, once again without loss of generality, that $\chi(x_0) = 1$, and  taking into account \cite[Theorem 6.1]{H}, we require further that $x_0$ satisfies \eqref{CdensFibulk}, \eqref{CdensFi}, where $\mathcal F_i(u;\cdot)$ is as in \eqref{auxf}.
	
	Choosing the sequence $\e_k \to 0^+$ in such a way that 
	$\mu(\partial Q(x_0,\e_k)) = 0$ and $Q(x_0,\e_k) \subset \Omega$,
	let $\{u_n\} \subset W^{2,1}(Q(x_0,\e_k);\Rb^N)$ be such that $u_n \to u$ in $W^{1,1}(Q(x_0,\e_k);\Rb^N)$ and
	\begin{equation}\label{F1Cdens}
		\dfrac{d\mathcal F_1(u; \cdot)}{d \lambda_u}(x_0) =
		\lim_{k \to +\infty}\frac{\mathcal F_1(u;Q(x_0,\e_k))}{\lambda_u(Q(x_0,\e_k))} 
		= \lim_{k,n}\frac{1}{\lambda_u(Q(x_0,\e_k))}
		\int_{Q(x_0,\e_k)}f_1(\mathcal \nabla^2 u_n(x)) \, dx.
	\end{equation}
	Then, as the constant sequence $\chi_n = \chi$ is admissible for 
	$\mathcal F(\chi,u;Q(x_0,\e_k))$, from \eqref{Dchi} and \eqref{F1Cdens} it follows that
	\begin{align}\label{estub}
		\frac{d \mu}{d \lambda_u}(x_0) &= \lim_{k \to +\infty}
		\frac{\mathcal F(\chi,u;Q(x_0,\e_k))}{\lambda_u(Q(x_0,\e_k))} \nonumber\\
		&\leq \liminf_{k,n}\left[\frac{1}{\lambda_u(Q(x_0,\e_k))}\int_{Q(x_0,\e_k)}
		f(\chi(x),\nabla^2 u_n(x))\,dx  +  |D \chi|(Q(x_0,\e_k))\right] \nonumber\\
		&\leq \lim_{k,n}\frac{1}{\lambda_u(Q(x_0,\e_k))}\int_{Q(x_0,\e_k)}
		f(1,\nabla^2 u_n(x))\,dx  \nonumber\\
		&\hspace{3cm} + \limsup_{k,n}\frac{1}{\lambda_u(Q(x_0,\e_k))}\int_{Q(x_0,\e_k)}
		\big[f(\chi(x),\nabla^2 u_n(x)) - f(1,\nabla^2 u_n(x))\big] \,dx\nonumber\\
		& = \dfrac{d\mathcal F_1(u; \cdot)}{d \lambda_u}(x_0) +
		\limsup_{k,n}\frac{1}{\lambda_u(Q(x_0,\e_k))}\int_{Q(x_0,\e_k)}
				\big[f(\chi(x),\nabla^2 u_n(x)) - f(1,\nabla^2 u_n(x))\big] \,dx.
	\end{align}
	The same argument used in the proof of Proposition~\ref{lbCantor}, namely a diagonal procedure together with Chacon's Biting Lemma,
yields
	$$\limsup_{k,n}\frac{1}{\lambda_u (Q(x_0,\e_k))}\int_{Q(x_0,\e_k)}
	\big[f(\chi(x),\nabla^2 u_n(x)) - f(1,\nabla^2 u_n(x))\big] \,dx = 0.$$
Hence, the conclusion follows from \eqref{estub}
by either \eqref{CdensFibulk} or \eqref{CdensFi} with $i=1$,
since
$$	\dfrac{d\mathcal F_1(u; \cdot)}{d \lambda_u}(x_0) = 	
\dfrac{d\mathcal F_1(u; \cdot)}{d \mathcal L^N}(x_0) = f(\chi(x_0),\nabla^2 u(x_0)),
$$
for $\mathcal L^N$ a.e. $x_0\in \Omega$,
or 
$$	\dfrac{d\mathcal F_1(u; \cdot)}{d \lambda_u}(x_0) = 	
\dfrac{d\mathcal F_1(u; \cdot)}{d |D^c(\nabla u)|}(x_0) = f^\infty\left(\chi(x_0),\frac{d D^c(\nabla u) }{d |D^c(\nabla u)|}(x_0)\right),
$$
for $|D^c(\nabla u)|$ a.e. $x_0 \in \Omega$.
\end{proof}

\subsection{The jump term}\label{jt}

Recall that, 
for every 
$(\chi, u)\in BV(\Omega;\{0,1\})\times BH(\Omega;\mathbb R^d)$, 
$\mathcal F(\chi,u,;\cdot)$ in \eqref{calFint} is the restriction to $\mathcal O(\Omega)$ of a Radon measure that we are denoting by $\mu$ and which 
can be decomposed as 
$$
\mu=\mu^a\mathcal L^N+ \mu^c+ \mu^j,
$$
where $\mu^j \ll |D^j(\nabla u)|+ |D \chi|$ and $\mu^c \ll |D^c(\nabla u)|$.

Our aim in this subsection is to characterise the measure $\mu^j$ and to this end we will follow the ideas of the global method for relaxation introduced in \cite{BFM}.

As before, we will continue to assume, without loss of generality, that $f$ is $2$-quasiconvex with respect to the second variable (cf. Proposition \ref{Frel=Frel**}).

Given $\chi \in BV(\O;\{0,1\})$, $u \in BH(\O;\mathbb R^d)$ and $A \in \mathcal O_{\infty}(\O)$,
we begin by defining the local Dirichlet-type functional
\begin{align}\label{mbis}
	m(\chi,u;A) := &\inf\Big\{\mathcal F (\theta,v;A) : 
	\theta \in  BV(\O;\{0,1\}), v \in BH(\O;\mathbb R^d),
	\theta\equiv \chi \hbox{ on }\partial A, \, {\rm supp}(v - u) \Subset A\Big\}. 
\end{align}

Our first objective, proved in Lemma \ref{BFT312} below, is to show that 
$\mathcal F(\chi,u,;A)$ and $m(\chi,u;A)$ behave in a similar fashion whenever $A$ is a cube of small side-length.

Fixing $\chi \in BV(\O;\{0,1\})$ and $u \in BH(\O;\mathbb R^d)$  
we define 
$$\lambda := \mathcal L^N + |(D^2)^s u| + |D\chi|$$ 
and, following \cite{BFM}, we let
$$\mathcal O^*:= \left\{Q_\nu(x,\e) : x \in \O,\,  \nu \in \mathbb S^{N-1},
\, \e > 0\right\}$$
and, for $\delta > 0$ and $A \in \mathcal O(\O)$, we set
\begin{align*}
	m^\delta(\chi,u;A) &:= \inf\Big\{\sum_{i = 1}^{+\infty}m(\chi,u;Q_i) : 
	Q_i \in \mathcal O^*, Q_i \cap Q_j = \emptyset \; {\rm if} \; i \neq j,\\
	& \hspace{4cm} Q_i \subset A, {\rm diam} \, Q_i < \delta, 
	\lambda\left(A \setminus \displaystyle\bigcup_{i=1}^{+\infty}Q_i\right) = 0\Big\}.
\end{align*}
Clearly, $\delta \mapsto m^\delta(\chi,u;A)$ is a decreasing function, so we define
$$m^*(\chi,u;A) := \sup \left\{m^\delta(\chi,u;A) : \delta > 0\right\}
= \lim_{\delta \to 0^+}m^\delta(\chi,u;A).$$

\begin{Lemma}\label{BFT311}
	Let $f$ be given by \eqref{densityint}, where $W_0$ and $W_1$ are continuous functions satisfying \eqref{growthint}.
	Given  $\chi \in BV(\O;\{0,1\})$, $u \in BH(\O;\mathbb R^d)$, we have
	$$\mathcal F(\chi,u;A) = m^*(\chi,u;A), \; \mbox{ for every } A \in \mathcal O(\O).$$
\end{Lemma}
\begin{proof}
The inequality 
$$m^*(\chi,u;A) \leq \mathcal F(\chi,u;A)$$ is standard. 
	
The proof of the reverse inequality relies on the lower semi-continuity of 
$\mathcal F(\cdot,\cdot;A)$ obtained in 
Proposition \ref{firstprop} iii) and on the definitions of $m^\delta(\chi,u;A)$, $m(\chi,u;A)$ and $m^*(\chi,u;A)$. Indeed, fixing $\delta > 0$, we consider $(Q_i^\delta)$ an admissible family for $m^\delta(\chi,u;A)$ such that, letting 
$\displaystyle N^\delta := A \setminus \displaystyle\cup_{i=1}^{+\infty}Q_i^\delta$,
\begin{equation}\label{est1}
\sum_{i =1}^{+\infty}m(\chi,u;Q_i^\delta) < m^\delta(\chi,u;A) 
+ \delta \; \mbox{ and } \; \lambda(N^\delta) = 0,
\end{equation}
and we now let $\theta_i^\delta \in BV(\O;\{0,1\})$ and 
$v_i^\delta \in BH(\O;\mathbb R^d)$
be such that $\theta_i^\delta = \chi$ on $\partial Q_i^\delta$, 
${\rm supp}(v_i^\delta - u) \Subset Q_i^\delta$ and
\begin{equation}\label{est2}
\mathcal F(\theta_i^\delta,v_i^\delta;Q_i^\delta) \leq 
m(\chi,u;Q_i^\delta) + \delta \mathcal L^N(Q_i^\delta).
\end{equation}
Setting 
$\displaystyle N_0^\delta := \O\setminus\displaystyle\cup_{i=1}^{+\infty}Q_i^\delta$,
we define 
$$\theta^\delta := \sum_{i =1}^{+\infty}\theta_i^\delta \,\chi_{Q_i^\delta} 
+ \chi \,\chi_{N_0^\delta} \;\;\; \mbox{ and } \; \;\;
v^\delta := \sum_{i =1}^{+\infty}v_i^\delta \,\chi_{Q_i^\delta} 
+ u \,\chi_{N_0^\delta}.$$
Following the computations in the proof of \cite[Lemma 3.11]{BFT} and
\cite[Lemma 4.2]{FHP}, we may show that $\theta^\delta \in BV(\O;\{0,1\})$, $v^\delta \in BH(\O;\mathbb R^d)$, 
$\theta^\delta \to \chi$  in $L^1(A;\{0,1\})$ and
$v^\delta \to u$ in $W^{1,1}(A;\Rb^d)$, as $\delta \to 0^+$.
	
Indeed, notice first that $v_i^\delta \in W^{1,1}(\O;\mathbb R^d)$, 
$\forall i \geq 1$, so the fact that ${\rm supp}(v_i^\delta - u) \Subset Q_i^\delta$ entails that $v^\delta \in W^{1,1}(\O;\mathbb R^d)$. On the other hand, an integration by parts argument proves that
$$D^2v^{\delta} = \sum_{i =1}^{+\infty}D^2v_i^{\delta} \lfloor Q_i^{\delta} +
D^2u \lfloor N_0^{\delta}.$$
Hence, the coercivity condition in Proposition \ref{firstprop} i), 
applied to each cube $Q_i^{\delta}$, together with \eqref{est1} and \eqref{est2}, ensure that $v^\delta \in BH(\O;\mathbb R^d)$. 
Using now the growth condition from above obtained in Proposition \ref{firstprop} i), we also conclude that
$$\mathcal F(\theta^\delta,v^\delta;N^\delta) \leq C \lambda(N^\delta) = 0.$$ 
Therefore, by the additivity of $\mathcal F(\theta^\delta,v^\delta;\cdot)$, we have
\begin{align}\label{est3}
\mathcal F(\theta^\delta,v^\delta;A) &= \sum_{i =1}^{+\infty}
\mathcal F(\theta_i^\delta,v_i^\delta;Q_i^\delta) + 
\mathcal F(\theta^\delta,v^\delta;N^\delta) \nonumber\\
&\leq \sum_{i =1}^{+\infty}m(\chi,u;Q_i^\delta) + \delta \mathcal L^N(A)
\leq m^\delta(\chi,u;A) + \delta + \delta \mathcal L^N(A).
\end{align}
	
To show that $v^\delta \to u$ in $W^{1,1}(A;\Rb^d)$, as $\delta \to 0^+$, we apply Poincaré's inequality to $\nabla u - \nabla v^{\delta}$ to obtain
\begin{align*}
\|\nabla u - \nabla  v^{\delta}\|_{L^1(\O;\mathbb R^{d\times N})} &=
\sum_{i =1}^{+\infty}\|\nabla u - \nabla  v_i^{\delta}\|_{L^1(Q_i^{\delta};\mathbb R^{d\times N})}\\
&\leq \sum_{i =1}^{+\infty}C\delta |D^2u - D^2v_i^{\delta}|(Q_i^{\delta})
\leq C\delta \big(|D^2u|(A) + |D^2v^{\delta}|(A)\big).
\end{align*}	
Again by coercivity, and \eqref{est3}, it follows that $ |D^2v^{\delta}|(A)$
is uniformly bounded in $\delta$ so we conclude that 
$\nabla v^{\delta} \to \nabla u$ in $L^1(\O;\mathbb R^{d \times N})$, 
as $\delta \to 0^+$.
Likewise, applying  Poincaré's inequality directly to $u - v^{\delta}$ yields
\begin{align*}
\|u - v^{\delta}\|_{L^1(\O;\mathbb R^{d})} &=
\sum_{i =1}^{+\infty}\|u - v_i^{\delta}\|_{L^1(Q_i^{\delta};\mathbb R^{d})}\\
&\leq \sum_{i =1}^{+\infty}C 
\|\nabla u - \nabla v_i^{\delta}\|_{L^1(Q_i^{\delta};\mathbb R^{d \times N})}
= C \|\nabla u - \nabla v^{\delta}\|_{L^1(\O;\mathbb R^{d \times N})},
\end{align*}
from where we conclude that 
$\displaystyle \lim_{\delta \to 0^+}\|u - v^{\delta}\|_{W^{1,1}(\O;\mathbb R^{d})} = 0$.
	
Given that $\theta^\delta \to \chi$ in $L^1(A;\{0,1\})$ and $v^\delta \to u$ in $W^{1,1}(A;\Rb^d)$, as $\delta \to 0^+$, the lower semi-continuity of 
$\mathcal F(\cdot,\cdot;A)$ allows us to conclude that 
\begin{align*}
\mathcal F(\chi,u;A) &\leq \liminf_{\delta \to 0^+}
\mathcal F(\theta^\delta,v^\delta;A) \\
&\leq \liminf_{\delta \to 0^+}
\left( m^\delta(\chi,u;A) + \delta + \delta \mathcal L^N(A)\right) = m^*(\chi,u;A)
\end{align*}
and this completes the proof.
\end{proof}

For $\delta > 0$, and given $A \in \mathcal O_{\infty}(\Omega)$, let 
$A_\delta:=\{x \in A: {\rm dist}(x, \partial A)> \delta\}$. Then, following
\cite[Lemma 4.1]{FHP} and \cite[Remark 3.2(1)]{BFM}, and taking also into account 
Remark \ref{A1A2} and Theorem \ref{thm5.7Hageryy}, 
 we conclude that
$$\limsup_{\delta \to 0^+}m(\chi,u,;A_\delta) \leq m(\chi,u;A).$$
Based on this fact, the previous lemma leads to the following
result whose proof is a straightforward adaptation of \cite[Lemma 3.5]{BFM}, see also \cite[Theorem 4.3]{FHP}.

\begin{Lemma}\label{BFT312}
Let $f$ be given by \eqref{densityint}, where $W_0$ and $W_1$ are continuous functions satisfying \eqref{growthint}. Let $ \mathcal F$ be as in \eqref{calFint} and $m$ as in \eqref{mbis}. 
Given  $\chi \in BV(\O;\{0,1\})$ and $u \in BH(\O;\mathbb R^d)$, we have
$$\lim_{\e \to 0^+}\frac{\mathcal F(\chi,u;Q_\nu(x_0,\e))}{\lambda(Q_\nu(x_0,\e))} = \lim_{\e \to 0^+}\frac{m(\chi,u;Q_\nu(x_0,\e))}{\lambda(Q_\nu(x_0,\e))},$$
for $\lambda$ a.e. $x_0 \in \O$ and for every $\nu \in \mathbb S^{N-1}$, where $\lambda:= \mathcal L^N+ |D \chi|+ |D^s(\nabla u)|$.
\end{Lemma}

The proof of the following result closely resembles the arguments in \cite[Proposition 4.5]{FHP}, with some modifications due to the presence of the $BV$ field $\chi$ in our functionals.

\begin{Lemma}\label{lemmattempt}
Let $\lambda$ be a bounded non-negative Radon measure on $\Omega$ and $m$ be as in \eqref{mbis}, with $F$ and $\mathcal F$ given in \eqref{Fint} and \eqref{calFint}, respectively, and  $W_i: \mathbb R_s^{d\times N \times N}\to \mathbb R, i=0,1$,  continuous functions satisfying \eqref{growthint}.
For $\varepsilon > 0$ and
given $u, v_\varepsilon \in BH(\Omega;\mathbb R^d)$ and 
$\chi, \eta_\varepsilon \in BV(\Omega;\{0,1\})$, 
let $x_0 \in \Omega$ and $\nu \in \mathbb S^{N-1}$  be such that
$\displaystyle	\lim_{\varepsilon \to 0^+}
\frac{m(\chi, u; Q_\nu (x_0, \varepsilon))}{\lambda(Q_\nu (x_0, \varepsilon))}$
exists. Then,
\begin{align*}
\hspace{-3cm}\limsup_{\varepsilon \to 0^+}
\frac{m(\eta_\varepsilon,v_\varepsilon ; Q_\nu (x_0, \varepsilon))}
{\lambda(Q_\nu (x_0, \varepsilon))}
-\lim_{\varepsilon \to 0^+}
\frac{m(\chi,u; Q_\nu (x_0, \varepsilon))}{\lambda(Q_\nu (x_0, \varepsilon))} 	
\end{align*}
\begin{align}\label{estimate}
&\hspace{1cm}\leq \limsup_{\delta \to 1^-}\limsup_{\varepsilon \to 0^+}
\frac{C}{\lambda(Q_\nu (x_0, \varepsilon))}
\Bigg\{\varepsilon^{N+1} + \varepsilon^N (1-\delta^N ) + 
|D^2u|(Q_\nu (x_0, \varepsilon)\setminus \overline{Q_\nu (x_0, \delta \varepsilon)}) \nonumber\\
&\hspace{2cm}+ |D^2v_\varepsilon|(Q_\nu (x_0, \varepsilon)\setminus 
\overline{Q_\nu (x_0, \delta \varepsilon)}) \nonumber \\
&\hspace{2cm}+\frac{1}{\varepsilon^2(1- \delta)^2}\int_{Q_\nu (x_0,\varepsilon)}
|u(x) - v_\varepsilon(x)| \,dx
+ \frac{1}{\varepsilon (1-\delta)}\int_{Q_\nu (x_0,\varepsilon)}
|\nabla u(x)-\nabla v_\varepsilon(x)| \,dx \nonumber \\
&\hspace{2cm}+ 
|D \eta_\varepsilon|(Q_\nu (x_0,\varepsilon)\setminus \overline{Q_{\nu}(x_0,\delta \varepsilon)}) 
+\int_{\partial Q_{\nu}(x_0,\delta \varepsilon)} 
|{\rm tr}(\chi)(x)-{\rm tr}(\eta_\varepsilon)(x)| \,d \mathcal H^{N-1}(x) \Bigg\}. 
\end{align}
\end{Lemma}
\begin{proof}
Fix $\delta \in (0, 1)$ and let $\varepsilon > 0$ be small enough so that 
$Q_\nu (x_0, \varepsilon) \subset \Omega$. Choose a cut-off function 
$\varphi \in C^\infty_c(Q_\nu(x_0, \varepsilon))$ such that $\varphi = 1$ in a neighbourhood of $Q_\nu (x_0, \varepsilon \delta)$ and
$$\displaystyle \|\nabla \varphi\|_{L^\infty}\leq \frac{ 2}{\varepsilon (1-\delta)},
\quad
\|\nabla ^2 \varphi\|_{L^\infty}\leq \frac{ 4}{\varepsilon^2 (1-\delta)^2}.$$
Now define
$$w_\varepsilon := \begin{cases}
\varphi \,u + (1 -\varphi)\,v_\varepsilon,  &\hbox{ in } Q_\nu (x_0, \varepsilon) \\
v_\varepsilon, &\hbox{ otherwise,}
\end{cases}$$
and choose 
$(\widehat \chi, \widehat u) \in BV(\Omega;\{0,1\})\times 
BH(\Omega;\mathbb R^d)$ such that
$\widehat \chi= \chi$ on $\partial Q_\nu (x_0, \varepsilon \delta)$, 
${\rm supp}(\widehat u-u) \Subset Q_\nu (x_0, \varepsilon \delta)$ and 
$$\varepsilon^{N+1} + m(\chi,u; Q_\nu (x_0, \varepsilon \delta)) 
\geq \mathcal F(\widehat \chi, \widehat u; Q_\nu (x_0, \varepsilon \delta)).$$
Set 
$$\widehat{v}_\varepsilon:=\begin{cases}
\widehat u, &\hbox{ in } Q_\nu(x_0,\varepsilon \delta)\\
w_\varepsilon &\hbox{ in }\Omega \setminus Q_\nu(x_0,\varepsilon \delta)
\end{cases} \qquad \mbox{and} \qquad
\widehat{\eta}_\varepsilon:=\begin{cases}
\widehat \chi, &\hbox{ in } Q_\nu(x_0, \varepsilon \delta)\\
\eta_\varepsilon, &\hbox{ in } \Omega \setminus Q_\nu(x_0,\varepsilon \delta).
\end{cases}$$
As $\widehat \eta_\varepsilon, \widehat v_\varepsilon$ are admissible for  $m(\eta_\varepsilon,v_\varepsilon ;Q_\nu(x_0,\varepsilon))$, by the locality of $\mathcal F$ and Remark \ref{A1A2}, we have
\begin{align}\label{est0}
m(\eta_\varepsilon, v_\varepsilon; Q_\nu(x_0, \varepsilon)) 
&\leq \mathcal F(\widehat \eta_\varepsilon, \widehat v_\varepsilon; 
Q_\nu (x_0, \varepsilon))\nonumber\\
&\leq \mathcal F(\widehat \chi, \widehat u; Q_\nu (x_0, \varepsilon \delta)) + 
C\mathcal L^N(Q_\nu (x_0, \varepsilon)\setminus Q_\nu (x_0, \varepsilon \delta))  \nonumber \\ 
&+ C \, |D^2 \widehat v_\varepsilon|
(Q_\nu (x_0, \varepsilon)\setminus Q_\nu (x_0, \varepsilon \delta)) 
+C \, |D \widehat \eta_\varepsilon|
(Q_\nu (x_0, \varepsilon)\setminus Q_\nu (x_0, \varepsilon \delta)) \nonumber\\
&\leq \varepsilon^{N+1} + m(\chi,u; Q_\nu (x_0, \varepsilon \delta))
+C\varepsilon^N (1- \delta^N ) \nonumber\\
&+ C\, |D^2w_\varepsilon|(Q_\nu (x_0, \varepsilon)\setminus \overline{Q_\nu (x_0, \varepsilon \delta)}) 
+ C \, |D\eta_\varepsilon|(Q_\nu(x_0,\varepsilon)\setminus
\overline{Q_{\nu}(x_0,\varepsilon \delta)})\nonumber\\ 
&+ C\int_{\partial Q_{\nu}(x_0, \varepsilon \delta)}
\Big(|{\rm tr}(\nabla \widehat u)(x)- {\rm tr}(\nabla w_\varepsilon)(x)|
+|{\rm tr}\,\widehat \chi(x)- {\rm tr}\,\eta_\varepsilon(x)|\Big) 
\,d \mathcal H^{N-1}(x). 
\end{align} 
		
Since, in $Q_\nu(x_0,\e)$,
$\nabla w_\varepsilon = (u - v_\varepsilon)\otimes \nabla \varphi 
+ \varphi\nabla u + (1 - \varphi)\nabla v_\varepsilon$,
we obtain
\begin{align*}
&|D^2w_\varepsilon|(Q_\nu(x_0,\varepsilon)\setminus 
\overline{Q_\nu (x_0, \varepsilon \delta)}) 
\leq C\bigg(|D^2u|(Q_\nu (x_0, \varepsilon)\setminus 
\overline{Q_\nu (x_0, \varepsilon \delta)})
+ |D^2v_\varepsilon|(Q_\nu (x_0, \varepsilon)\setminus 
\overline{Q_\nu (x_0, \varepsilon \delta)})\\
& \hspace{3,8cm}+ \frac{1}{\varepsilon^2(1 -\delta)^2}
\int_{Q_\nu (x_0,\varepsilon)}|u(x) - v_\varepsilon(x)| \, dx
+ \frac{1}{\varepsilon(1 - \delta)}
\int_{Q_\nu (x_0,\varepsilon)}|\nabla u(x)-\nabla v_\varepsilon(x)|\, dx\bigg).
\end{align*}
Observe also that ${\rm tr}(\nabla \widehat u)= {\rm tr}(\nabla u)= {\rm tr}(\nabla w_\varepsilon)$ on $\partial Q_{\nu}(x_0, \varepsilon \delta)$.
	
Furthermore, from Lemma \ref{FM}, we deduce that
\begin{align*}
\limsup_{\delta \to 1^{-}}\limsup_{\varepsilon \to 0^+}
\frac{m(\chi, u; Q_\nu (x_0, \varepsilon \delta))}
{\lambda(Q_\nu (x_0, \varepsilon))} &= 
\limsup_{\delta \to 1^{-}}\limsup_{\varepsilon \to 0^+}
\left[\frac{m(\chi, u; Q_\nu (x_0, \varepsilon \delta))}
{\lambda(Q_\nu (x_0, \varepsilon \delta))} 
\frac{\lambda(Q_\nu (x_0, \varepsilon \delta))}
{\lambda(Q_\nu (x_0, \varepsilon))}\right]
\\
&\leq \lim_{\varepsilon \to 0^+}\frac{m(\chi, u; Q_\nu(x_0, \varepsilon))}
{\lambda(Q_\nu (x_0, \varepsilon))},
\end{align*}
hence, to complete the proof, it suffices to divide inequality \eqref{est0}
by $\lambda(Q_\nu (x_0, \varepsilon))$ and take the $\limsup$ as $\varepsilon \to 0^+$  and $\delta \to 1^-$, using the above estimate in the appropriate term.
\end{proof}

For $(a,b,C,D,\nu)  \in \{0,1\} \times \{0,1\} \times \mathbb{R}^{d\times N} \times\mathbb{R}^{d\times N} \times \mathbb S^{N-1}$, with  $C-D= P \otimes \nu$, for some $P \in \mathbb R^d$, for $y \in \Omega$, we define the functions $\chi_{a,b,\nu}$ and $u_{C,D,\nu}$ by
\begin{equation}\label{bounddata}
\chi_{a,b,\nu}(y) := \begin{cases}
a, & {\rm if } \; y \cdot \nu \geq 0\\
b, & {\rm if } \; y \cdot \nu < 0
\end{cases}
\; \; \; {\rm and } \; \; \;
u_{C,D,\nu}(y) := \begin{cases}
C \cdot y, & {\rm if } \; y \cdot \nu \geq 0\\
D \cdot y, & {\rm if } \; y \cdot \nu < 0.
\end{cases}
\end{equation}

Given  $\chi \in BV(\Omega;\{0,1\})$ and $u \in BH(\Omega;\mathbb R^d)$, for 
$x_0 \in S_{\chi} \cup S_{\nabla u}$, $\nu(x_0)$ will be used to denote $\nu_{\chi}(x_0)$, if $x_0 \in S_{\chi} \setminus S_{\nabla u}$, and 
$\nu_{\nabla u}(x_0)$, if $x_0 \in S_{\nabla u} \setminus S_{\chi}$. On the other hand, since  $(\chi, \nabla u)$ may be considered as components of the same $SBV$ vector-valued field, for $\mathcal H^{N-1}$ a.e. $x_0 \in S_\chi \cap S_{\nabla u}$ we may select $\nu(x_0):= \nu_{\chi}(x_0)=\nu_{\nabla u}(x_0)$ where the orientation of $\nu_{\chi}(x_0)$ is chosen so that $\chi^+(x_0)=1, \chi^{-}(x_0)=0$ and then $(\nabla u)^+(x_0), (\nabla u)^-(x_0)$ are selected according to this orientation,
see \cite{AFP}. This is the definition of $\nu(x_0)$ that we will consider in what follows.

\begin{Proposition}\label{mg}
Let $f$ be given by \eqref{densityint}, where $W_0$ and $W_1$ are continuous functions satisfying \eqref{growthint}, and $m$ be as in \eqref{mbis}, with $F$ and $\mathcal F$ given in \eqref{Fint} and \eqref{calFint}, respectively.
Given $\chi \in BV(\Omega;\{0,1\})$ and $u \in BH(\Omega;\mathbb R^d)$,  we have
$$\mathcal{F}( \chi, u; \Omega \cap (S_\chi \cup S_{\nabla u})) =
\int_{\Omega \cap (S_\chi \cup S_{\nabla u})} 
g(x,\chi^+(x),\chi^-(x),u(x),(\nabla u)^+(x), (\nabla u)^-(x),\nu(x))
\, d\mathcal H^{N-1}(x),$$
where $\nu(x)$ denotes the unit normal to $S_\chi \cup S_{\nabla u}$,
\begin{equation}\label{defg}
g(x_0,a,b,r, C,D,\nu) := \limsup_{\varepsilon \to 0^+}
\frac{m(\chi_{a,b,\nu}(\cdot - x_0),r+ u_{C,D,\nu}(\cdot - x_0);Q_\nu(x_0,\e))}
{\e^{N-1}}
\end{equation}
and $\chi_{a,b,\nu}$, $u_{C,D,\nu}$ are the functions defined in \eqref{bounddata}, for $a,b \in \{0,1\}, r \in \mathbb R^d$, and $C,D \in \mathbb R^{d \times N},$
with  $C-D= P \otimes \nu$, for some $P \in \mathbb R^d$.
\end{Proposition}
\begin{proof}
Given $x_0 \in S_\chi \cup S_{\nabla u}$, in the sequel, for simplicity of notation, we will write $\nu = \nu(x_0)$, and without loss of generality we will assume that 
$u \in SBH(\Omega;\mathbb R^d)$. 

Let $x_0 \in \O \cap (S_\chi \cup S_{\nabla u})$ be a point satisfying the following properties
\begin{align}
&\lim_{\varepsilon \rightarrow 0^+}\frac{1}{\e^N}
\int_{Q_{\nu}(x_0,\e)}|\chi(x) - {\chi}(x_0)| \, dx = 0, 
\; {\rm if} \;x_0 \in \Omega \setminus S_{\chi},\nonumber\\
&\lim_{\varepsilon \rightarrow 0^+}\frac{1}{\e^N}
\int_{Q^+_{\nu}(x_0,\e)}\hspace{-0,2cm}|\chi(x) - {\chi}^+(x_0)| \, dx = 
\lim_{\varepsilon \rightarrow 0^+}\frac{1}{\e^N}
\int_{Q^-_{\nu}(x_0,\e)}\hspace{-0,2cm}|\chi(x) - {\chi}^-(x_0)| \, dx =0, 
\; {\rm if} \;x_0 \in \Omega \cap S_{\chi}, \label{chipm}\\
&\lim_{\varepsilon \rightarrow 0^+}\frac{1}{\e^N}
\int_{Q_{\nu}(x_0,\e)}|\nabla u(x) - {\nabla u}(x_0)| \, dx = 0, 
\; {\rm if} \;x_0 \in \Omega \setminus S_{\nabla u}, \nonumber\\
&\lim_{\varepsilon \rightarrow 0^+}\frac{1}{\e^N}
\int_{Q^+_{\nu}(x_0,\e)}\hspace{-0,2cm}|\nabla u(x) - {(\nabla u)}^+(x_0)| \, dx = 
\lim_{\varepsilon \rightarrow 0^+}\frac{1}{\e^N}
\int_{Q^-_{\nu}(x_0,\e)}\hspace{-0,2cm}|\nabla u(x) - {(\nabla u)}^-(x_0)| \, dx =0, 
\; {\rm if} \;x_0 \in \Omega \cap S_{\nabla u},\label{upm}
\end{align}
where 
$$Q^{\pm}_{\nu}(x_0,\e) = 
\left\{x \in Q_{\nu}(x_0,\e) : (x - x_0) \cdot (\pm \nu) > 0\right\},$$
and, with an abuse of notation, we use the same symbol for $u$ and its approximate limit $\widetilde u$, for $\nabla u$ and its approximate limit 
$\widetilde{\nabla u}$, in $\Omega\setminus S_{\nabla u}$, and for $\chi$ and its approximate limit $\widetilde \chi$, in $\Omega \setminus S_\chi$ (see \cite[Definition 3.63 and Proposition 3.64]{AFP}).

Defining, for every $x \in Q_\nu(x_0,\varepsilon)$,
\begin{align*}
v_J(x):= u(x_0)+ 
\begin{cases}
(\nabla u)^+(x_0) \cdot (x-x_0) & \hbox{ if } (x-x_0)\cdot \nu \geq 0,\\
(\nabla u)^-(x_0) \cdot (x-x_0) & \hbox{ if } (x-x_0)\cdot \nu <0	
\end{cases}
\end{align*}
and
\begin{align}\label{chi0def}
\chi_J(x) := \begin{cases}
\chi^+(x_0), & {\rm if } \;  (x-x_0)\cdot \nu \geq 0,\\
\chi^-(x_0), & {\rm if } \; (x-x_0) \cdot \nu < 0
\end{cases}
\end{align}
then, \cite[Theorem 2.1]{FHP} entails that 
\begin{equation}\label{vJ}
\lim_{\varepsilon \rightarrow 0^+}\frac{1}{\e^{N+1}} 
\int_{Q_{\nu}(x_0,\e)}|u(x) - v_J(x)| \, dx = 0, 
\; {\rm if} \;x_0 \in \Omega \cap S_{\nabla u},
\end{equation}
whereas from \eqref{upm} we obtain 
\begin{equation}\label{nablavJ}
\lim_{\varepsilon \rightarrow 0^+}\frac{1}{\e^{N}}
\int_{Q_{\nu}(x_0,\e)}|\nabla u(x) - {\nabla v_J}(x)| \, dx = 0, 
\; {\rm if} \;x_0 \in \Omega \cap S_{\nabla u}.
\end{equation}
Furthermore, the point $x_0$ is also chosen so that
\begin{equation}\label{symjump}
\lim_{\varepsilon \to 0}\frac{|D^2 u|(Q_\nu(x_0,\varepsilon))}{\varepsilon^{N-1}}
=|[\nabla u](x_0)| =|[\nabla v_J](x_0)|,  \hbox{ if } x_0 \in S_{\nabla u},	\end{equation}
\begin{equation}\label{perchi0}
\lim_{\e \to 0^+}\frac{1}{\e^{N-1}}|D\chi|(Q_{\nu}(x_0,\e)) = |[\chi](x_0)|= 1,
\end{equation}
and	
\begin{equation}\label{muj}
\mu^j(x_0) = \lim_{\varepsilon \rightarrow 0^+}
\frac{\mathcal F(\chi,u;Q_{\nu}(x_0,\e))}
{\mathcal H^{N-1} \lfloor (S_\chi \cup S_{\nabla u})(Q_{\nu}(x_0,\e))} =
\lim_{\varepsilon \rightarrow 0^+}\frac{1}{\e^{N-1}}\int_{Q_{\nu}(x_0,\e)}d\mu(x)
\mbox{ exists and is finite}.
\end{equation}
These properties hold for $\mathcal H^{N-1}$ a.e. 
$x_0 \in \Omega\cap (S_\chi \cup S_{\nabla u})$.

Notice also that, for $x \in Q_\nu(x_0,\e)$, $\chi_{a,b,\nu}(\cdot - x_0) = \chi_J$,
when $a= \chi^+(x_0), b = \chi^-(x_0)$ and, when $r= u(x_0)$, $C= (\nabla u)^+(x_0)$,
$D= (\nabla u)^-(x_0)$, we have $r + u_{C,D,\nu}(\cdot - x_0) = v_J$.

Letting $\sigma := \mathcal H^{N-1}\lfloor (S_\chi \cup S_{\nabla u})$, by Lemma \ref{BFT312} it follows that, for $\sigma$ a.e. $x_0 \in \O$,
\begin{equation}\label{msigma}
\frac{d\mathcal F(\chi,u;\cdot)}{d \sigma}(x_0) =
\lim_{\varepsilon \to 0^+}
\frac{\mathcal F(\chi,u;Q_\nu(x_0,\e))}{\sigma(Q_\nu(x_0,\e))} =
\lim_{\varepsilon \to 0^+}
\frac{m(\chi,u;Q_\nu(x_0,\e))}{\sigma(Q_\nu(x_0,\e))}
\end{equation}
and
\begin{equation}\label{msigma2}
\lim_{\varepsilon \to 0^+}
\frac{\mathcal F(\chi_J,v_J;Q_\nu(x_0,\e))}{\sigma(Q_\nu(x_0,\e))} =
\lim_{\varepsilon \to 0^+}
\frac{m(\chi_J,v_J;Q_\nu(x_0,\e))}{\sigma(Q_\nu(x_0,\e))}.
\end{equation}

Hence, from \eqref{muj}-\eqref{msigma2}, we obtain
\begin{align*}
\frac{d\mathcal F(\chi,u;\cdot)}{d \sigma}(x_0) &= 
\lim_{\varepsilon \to 0^+}
\frac{m(\chi,u;Q_\nu(x_0,\e))}{\sigma(Q_\nu(x_0,\e))} \\
& = \lim_{\varepsilon \to 0^+}\frac{m(\chi,u;Q_\nu(x_0,\e))- 
m(\chi_J,v_J;Q_\nu(x_0,\e)) 
+ m(\chi_J,v_J;Q_\nu(x_0,\e))}{\e^{N-1}}\\
& = \lim_{\varepsilon \to 0^+}
\frac{m(\chi_J,v_J;Q_\nu(x_0,\e))}{\e^{N-1}} = 
g(x_0,\chi^+(x_0),\chi^-(x_0),u(x_0),(\nabla u)^+(x_0), (\nabla u)^-(x_0),\nu(x_0)),
\end{align*}
because
$$\lim_{\varepsilon \to 0^+}\frac{m(\chi,u;Q_\nu(x_0,\e))}{\e^{N-1}}- 
\lim_{\varepsilon \to 0^+}\frac{m(\chi_J,v_J;Q_\nu(x_0,\e))}{\e^{N-1}} =0.$$
Indeed,	we apply Lemma \ref{lemmattempt} to the measure $\lambda = \sigma$ and to the constant sequence $(\chi_J,v_J)$. Thus, from \eqref{estimate}, 
using the fact that 
$$\frac{|D^2 v_J|(Q_\nu (x_0, \varepsilon)\setminus \overline{Q_\nu (x_0, \delta \varepsilon)})}{\varepsilon^{N-1}}= \big|[\nabla u](x_0)\big|(1 -\delta^{N-1}),$$ 
and 
$$\frac{|D \chi_J|(Q_\nu(x_0,\varepsilon)\setminus \overline{Q_{\nu}(x_0,\delta \varepsilon)})}{\varepsilon^{N-1}} =\big |[\chi(x_0)]\big|(1-\delta^{N-1}), $$ 
we have
\begin{align*}
&\left|\lim_{\varepsilon \to 0^+}
\frac{m(\chi,u;Q_\nu(x_0,\varepsilon))}{\varepsilon^{N-1}}
- \lim_{\varepsilon \to 0^+}
\frac{m(\chi_J,v_J;Q_\nu(x_0,\varepsilon))}{\varepsilon^{N-1}}\right|\\
&\hspace{1cm}\leq C \limsup_{\delta \to 1^-}\limsup_{\varepsilon \to 0^+}
\bigg\{\varepsilon^2 + \varepsilon(1-\delta^N ) +
\frac{|D^2u|(Q_\nu (x_0, \varepsilon)\setminus 
\overline{Q_\nu (x_0, \delta \varepsilon)})}{\varepsilon^{N-1}} 
+\big|[\nabla u](x_0)\big|(1 -\delta^{N-1})\\
&\hspace{1cm} + \frac{1}{(1-\delta)^2}\frac{1}{\varepsilon^{N+1}}
\int_{Q_\nu(x_0,\varepsilon)}|u(x)- v_J(x)| \, dx
+ \frac{1}{(1-\delta)\varepsilon^N} 
\int_{Q_\nu(x_0,\varepsilon)}|\nabla u(x)-\nabla v_J (x)| \, dx\\
&\hspace{1cm} + \big|[\chi(x_0)]\big|(1-\delta^{N-1} ) 
+ \int_{\partial Q_{\nu}(x_0, \varepsilon \delta))}
|{\rm tr}(\chi)(x)- {\rm tr}(\chi_J)(x)| \,d {\mathcal H}^{N-1}(x) \bigg\}=0.
\end{align*}
This is due to \eqref{vJ}, \eqref{nablavJ}, taking into account that \eqref{symjump} guarantees that
\begin{align*}
0\leq  \limsup_{\delta \to 1^-}\limsup_{\varepsilon \to 0^+}
\frac{|D^2u|(Q_\nu (x_0, \varepsilon)\setminus \overline{Q_\nu (x_0, \delta \varepsilon)})}{\varepsilon^{N-1}}
\leq \limsup_{\delta \to 1^-}\big|[\nabla u](x_0)\big|(1- \delta^{N-1}) = 0,
\end{align*}
and, also, since
$$\limsup_{\delta \to 1^-}\limsup_{\varepsilon \to 0^+}
\int_{\partial Q_{\nu}(x_0, \varepsilon \delta)}
|{\rm tr}(\chi)(x)- {\rm tr}(\chi_J)(x)| \,d {\mathcal H}^{N-1}(x) =0,$$ 
which can be proved, arguing as in \cite[eq. (4.43)]{BMZ}, using
the continuity of the trace operator with respect to $BV$ strict convergence, which is ensured by \eqref{chipm}, \eqref{chi0def} and \eqref{perchi0}.

Therefore, we conclude that
\begin{align*}
\mathcal{F}( \chi, u; \Omega \cap (S_\chi \cup S_{\nabla u})) &=
\int_{\Omega \cap (S_\chi \cup S_{\nabla u})} 
\frac{d\mathcal F(\chi,u;\cdot)}{d \sigma}(x) \, d\sigma(x) \\
&=\int_{\Omega \cap (S_\chi \cup S_{\nabla u})} 
g(x,\chi^+(x),\chi^-(x),u(x),(\nabla u)^+(x),(\nabla u)^-(x),\nu(x))\, d\mathcal H^{N-1}(x),
\end{align*}
and the proof is complete.
\end{proof}

\begin{Remark}\label{remGM}{\rm 
We now make some comments on the previous results.
\begin{itemize}
\item[i)] The arguments presented in Lemma \ref{BFT312}, Lemma \ref{lemmattempt} and Proposition \ref{mg} can be applied to the case of a relaxed functional, as in \eqref{calFint}, based on the energy  in \eqref{enernew},  or even to generic functionals 
$\mathcal F:BV(\Omega;\{0,1\})\times BH(\Omega;\mathbb R^d)\times 
\mathcal O(\Omega) \to [0,+\infty]$ satisfying i), iii) and iv) of Proposition \ref{firstprop} and such that 
$\mathcal F(\chi, u; \cdot)$ is the restriction of a Radon measure.
In particular, this applies to $\mathcal F$ as in \eqref{calFint}, with
$F:BV(\Omega;\{0,1\})\times W^{2,1}(\Omega;\mathbb R^d)\times \mathcal O(\Omega) 
\to [0,+\infty]$  given by
\begin{align}\label{F2}
F(\chi,u;A)= \int_A f_1\color{black}(x,\chi(x), u(x), \nabla u(x), \nabla^2 u(x)) \, dx +|D\chi|(A),
\end{align}
where $f_1$ is continuous and satisfies linear growth and coercivity conditions in the last variable.
Moreover, arguing as in \cite[Section 6]{FHP}, one can show that for every 
$(\chi, u) \in BV(\Omega;\{0,1\})\times SBH(\Omega;\mathbb R^d)$ and every 
$A \in \mathcal O(\Omega)$, 
\begin{align*}
\mathcal F(\chi,u;A) &=\int_A h(x,\chi(x),u(x),\nabla u(x), \nabla^2 u(x)) \, dx \nonumber\\
&+ \int_{A \cap \left(S_{\chi}\cup S_{\nabla u}\right)}
g_1\color{black}(x,\chi^+(x), \chi^-(x), u(x), (\nabla u)^+(x),(\nabla u)^-(x),
\nu(x))\, d \mathcal H^{N-1}(x), 
\end{align*}
where 
\begin{align}\label{hdef}
h(x_0,\rho,r,\xi,\Sigma):=\limsup_{\varepsilon \to 0^+}\frac{m(\rho, r+\xi(\cdot-x_0)+\tfrac{1}{2}\Sigma (\cdot-x_0, \cdot-x_0);Q(x_0,\varepsilon))}{\varepsilon^N},
\end{align} 
for all $x_0\in \Omega, \rho \in \{0,1\}, r \in \mathbb R^d, \xi \in \mathbb R^{d\times N}, \Sigma \in \mathbb R^{d\times N \times N}_s$. Here
$g_1$ is as in \eqref{defg}, with $m$ defined in \eqref{mbis}, associated to 
$\mathcal F$ given by \eqref{calFint} for the energy $F$ in \eqref{F2}.
\item[ii)] In the case under consideration, that is, when $F$ takes the form in \eqref{Fint}, it is easy to see that $g_1$ does not depend explicitly on either $x$ or $u(x)$ (see also \cite{FHP}). Indeed, it suffices to observe that, in this case, for all
$(x_0,a,b,r,C,D,\nu)  \in \Omega \times \{0,1\} \times \{0,1\}  \times \mathbb R^d \times \mathbb{R}^{d\times N} \times\mathbb{R}^{d\times N} \times \mathbb S^{N-1}$,
we have
$$g_1(x_0,a,b,r,C,D,\nu) = g_1(0,a,b,0,C,D,\nu),$$ 
 which we denote by $g(a,b,C,D,\nu)$. 

The proof of the previous equality amounts to showing that
$$m\big(\chi_{a,b,\nu}(\cdot -x_0),r + u_{C,D,\nu}(\cdot -x_0);Q_\nu(x_0,\varepsilon)\big) =
m\big(\chi_{a,b,\nu},u_{C,D,\nu};Q_\nu(0,\varepsilon)\big),$$ 
where $\chi_{a,b,\nu}$ and $u_{C,D,\nu}$ are given in \eqref{bounddata},
and this relies on a change of variables argument using translations. 
In fact, if a pair $(\theta,v)$ is admissible for
$m\big(\chi_{a,b,\nu}(\cdot -x_0),r + u_{C,D,\nu}(\cdot -x_0);Q_\nu(x_0,\varepsilon)\big)$, then the pair
$(\eta,w)$ defined by $\eta(y) = \theta(x_0 +y)$, $w(y) = v(x_0 +y) - r$, 
for $y \in Q_\nu(0,\varepsilon)$, is admissible for
$m\big(\chi_{a,b,\nu},u_{C,D,\nu};Q_\nu(0,\varepsilon)\big)$,
and similarly for the reverse case.
\item[iii)] It is also routine to verify that, if $\chi^+(x_0)=\chi^-(x_0)$, then (see \cite{H})
\begin{align*}
g(\chi^+(x_0),\chi^-(x_0),(\nabla u)^+(x_0), (\nabla u)^- (x_0),\nu(x_0)) &
=  Q^2(f^\infty)(\chi(x_0),[\nabla u](x_0)\otimes \nu(x_0)) \\
& =  (Q^2f)^\infty(\chi(x_0),[\nabla u](x_0)\otimes \nu(x_0)),
\end{align*}
where Proposition \ref{SQfinfty=} has been taken into account, and if $\nabla u^+(x_0)=\nabla u^-(x_0)$, then 
$$g(\chi^+(x_0),\chi^-(x_0),\nabla u^+(x_0), \nabla u^- (x_0),\nu(x_0))= [\chi](x_0).$$
\item[(iv)] Finally, it is worth observing that, when $f_1$ in \eqref{F2} has the form given by \eqref{enernew}, then, standard arguments show that 
$h(x, \chi(x),u(x),\nabla u(x),\nabla^2 u(x))$ in \eqref{hdef} coincides with $$Q^2f(\chi(x), \nabla^2 u(x)) + W_2(x,\chi(x),u(x), \nabla u(x)),$$ 
as in \eqref{fing1}.
\end{itemize}}
\end{Remark}

Propositions \ref{lbCantor}, \ref{ub}, \ref{mg} and part ii) of Remark \ref{remGM} 
provide the proof of our main result stated in Theorem \ref{main}.

Taking into account Remark \ref{remGM} ii), our next result provides a more explicit, sequential, characterisation of the relaxed surface energy density $g$, which is independent of the Dirichlet functional 
\eqref{mbis}. It is obtained under the assumption that \eqref{finfty} holds.

\begin{Proposition}\label{Ktilde}
Let $f$ be given by \eqref{densityint}, where $W_0$ and $W_1$ are continuous functions satisfying \eqref{growthint}, and assume that $f$ is $2$-quasiconvex in the second variable and that \eqref{finfty} holds.
For every $(a,b,C,D,\nu)  \in \{0,1\} \times \{0,1\}  \times \mathbb{R}^{d \times N} \times\mathbb{R}^{d \times N} \times \mathbb S^{N-1}$, such that $C-D=P\otimes \nu$, for some $P\in\mathbb R^d$, we have
$$g(a,b,C,D,\nu) = \widetilde{K}(a,b,C,D,\nu),$$
where $g$ is as in \eqref{defg},
\begin{align*}
\widetilde{K}(a,b,C,D,\nu) & := \inf\bigg\{\liminf_{n \to + \infty}
\left[\displaystyle\int_{Q_{\nu}}f^{\infty}(\chi_n(x),\nabla^2 u_n(x)) \, dx
+ |D\chi_n|(Q_{\nu})\right] : \{\chi_n\} \subset BV\left(Q_{\nu};\{0,1\}\right), \nonumber \\
&\hspace{1,1cm} \{u_n\} \subset W^{2,1}\left(Q_{\nu};\mathbb{R}^d\right), 
\chi_n \to \chi_{a,b,\nu}
\; {\rm in } \; L^1(Q_{\nu};\{0,1\}), u_n \to u_{C,D,\nu} \; {\rm in } \; W^{1,1}(Q_{\nu};\mathbb R^d) \bigg\}, 
\end{align*}
and $\chi_{a,b,\nu}$, $u_{C,D,\nu}$ are defined in \eqref{bounddata}.
\end{Proposition}
\begin{proof}
By Theorem \ref{main}, applied to $\chi_{a,b,\nu}$ and $u_{C,D,\nu}$ in $Q_\nu(0,\e)$, we have
$$\mathcal F(\chi_{a,b,\nu},u_{C,D,\nu};Q_\nu(0,\e)) = f(a,0)\,\frac{\e^N}{2}
+ f(b,0)\,\frac{\e^N}{2} + g(a,b,C,D,\nu) \,\e^{N-1}$$
so that, from \eqref{calFint}, we obtain 
\begin{align*}
g(a,b,C,D,\nu) & = \lim_{\varepsilon \to 0^+}
\frac{\mathcal F(\chi_{a,b,\nu},u_{C,D,\nu};Q_{\nu}(0,\varepsilon))}{\e^{N-1}}\nonumber \\
& = \lim_{\varepsilon \to 0^+}\frac{1}{\varepsilon^{N-1}}
\inf\bigg\{\liminf_{n \to + \infty}
\left[\displaystyle\int_{Q_{\nu}(0,\varepsilon)}f(\chi_n(x),\nabla^2 u_n(x)) \, dx
+ |D\chi_n|(Q_{\nu}(0,\varepsilon))\right] : \nonumber \\
& \hspace{3cm} \{\chi_n\} \subset BV\left(Q_{\nu}(0,\varepsilon);\{0,1\}\right), 
\{u_n\} \subset W^{2,1}\left(Q_{\nu}(0,\varepsilon);\mathbb{R}^d\right),  \nonumber \\
& \hspace{3cm}\chi_n \to \chi_{a,b,\nu}
\; {\rm in } \; L^1(Q_{\nu}(0,\varepsilon);\{0,1\}), u_n \to u_{C,D,\nu} \; {\rm in } \; W^{1,1}(Q_{\nu}(0,\varepsilon);\mathbb R^d) \bigg\} \nonumber \\
& = \lim_{\varepsilon \to 0^+}\frac{1}{\varepsilon^{N-1}}\lim_{n \to + \infty}
\left[\displaystyle\int_{Q_{\nu}(0,\varepsilon)}f(\chi_{n,\varepsilon}(x),
\nabla^2 u_{n,\varepsilon}(x)) \, dx
+ |D\chi_{n,\varepsilon}|(Q_{\nu}(0,\varepsilon))\right], 
\end{align*}
where, in view of Proposition \ref{firstprop} ii), $\{\chi_{n,\varepsilon}\} \subset BV\left(Q_{\nu}(0,\varepsilon);\{0,1\}\right)$ and $\{u_{n,\varepsilon}\} \subset W^{2,1}\left(Q_{\nu}(0,\varepsilon);\mathbb{R}^d\right)$ are recovery sequences for 
${\mathcal F}(\chi_{a,b,\nu},u_{C,D,\nu};Q_\nu(0,\varepsilon))$. In particular, for every $\e > 0$, $\chi_{n,\e} \to \chi_{a,b,\nu}$ in
$L^1(Q_{\nu}(0,\varepsilon);\{0,1\})$ and $u_n \to u_{C,D,\nu}$ 
in $W^{1,1}(Q_{\nu}(0,\varepsilon);\mathbb R^d)$, as $n \to + \infty$.
 
Defining, for $y \in Q_\nu$, $\theta_{n,\varepsilon}(y):= \chi_{n,\varepsilon}(\varepsilon y)$ and 
$\displaystyle v_{n,\varepsilon}(y)=\frac{1}{\varepsilon} u_{n,\varepsilon}(\varepsilon y)$, it follows that 
for every $\varepsilon$, $\theta_{n,\varepsilon} \to \chi_{a,b,\nu}$ in $L^1(Q_\nu;\{0,1\})$, as $n\to +\infty$, and $v_{n,\varepsilon}\to u_{C,D,\nu}$ in $W^{1,1}(Q_\nu, \mathbb R^d)$, as $n\to +\infty$, and a simple change of variables now yields
\begin{align}\label{secondg}
g(a,b,C,D,\nu) & = \lim_{\varepsilon \to 0^+}\lim_{n\to +\infty}
\left[\int_{Q_{\nu}}\varepsilon f\left(\theta_{n,\varepsilon}(y),
\frac{1}{\varepsilon}\nabla^2 v_{n,\varepsilon}(y)\right) \, dy
+ |D\theta_{n,\varepsilon}|(Q_{\nu})\right]\nonumber\\
&=\lim_{\varepsilon \to 0^+}\lim_{n\to +\infty}
\left[\displaystyle\int_{Q_{\nu}}
f^\infty(\theta_{n,\varepsilon}(y),\nabla^2 v_{n,\varepsilon}(y)) \, dy
+ |D\theta_{n,\varepsilon}|(Q_{\nu})\right]
\end{align}
since, by \eqref{finfty}, we have
\begin{equation}\label{limitffinfty}
\lim_{\varepsilon \to 0^+}\lim_{n\to +\infty}\int_{Q_\nu}\left[
\varepsilon f\left(\theta_{n,\varepsilon}(y),
\frac{1}{\varepsilon}\nabla^2 v_{n,\varepsilon}(y)\right) - 
f^\infty(\theta_{n,\varepsilon}(y),\nabla^2 v_{n,\varepsilon}(y))\right] \, dy = 0.
\end{equation}

Indeed, we may write
\begin{align*}
& \int_{Q_{\nu}}\left[\e f\left(\theta_{n,\e}(y), \frac{1}{\e}\nabla^2 v_{n,\e}(y)\right)- f^{\infty}(\theta_{n,\e}(y),\nabla^2 v_{n,\e}(y)) \right]\, dy \\
& = \int_{Q_{\nu}\cap \{\frac{1}{\e}|\nabla^2 v_{n,\e}(y)| \leq L\}}
\left[\e f\left(\theta_{n,\e}(y), \frac{1}{\e}\nabla^2 v_{n,\e}(y)\right)- f^{\infty}(\theta_{n,\e}(y),\nabla^2 v_{n,\e}(y))\right] \, dy \\
& + \int_{Q_{\nu}\cap \{\frac{1}{\e}|\nabla^2 v_{n,\e}(y)| > L\}}
\left[\e f\left(\theta_{n,\e}(y), \frac{1}{\e}\nabla^2 v_{n,\e}(y)\right)- f^{\infty}(\theta_{n,\e}(y),\nabla^2 v_{n,\e}(y))\right] \, dy
= : I_1 + I_2.
\end{align*}
By \eqref{growthint}, \eqref{densityint} and the linear growth, in the second variable, of $f^\infty$, easily deduced from that of $f$, we have
\begin{align*}
|I_1| &\leq \int_{Q_{\nu}\cap \{|\nabla^2 v_{n,\e}(y)| \leq \e L\}}
\left[\e \, C \left(1 + \frac{1}{\e}|\nabla^2 v_{n,\e}(y)|\right) 
+ C \, |\nabla^2 v_{n,\e}(y)| \right]\, dy \\
& \leq \int_{Q_{\nu}}\e \, C \, dy = O(\e)
\end{align*}
and, by hypothesis \eqref{finfty} with $t = \frac{1}{\e}$, H\"older's inequality and
the growth of $f$, it follows that
\begin{align*}
|I_2| &\leq \int_{Q_{\nu}\cap \{\frac{1}{\e}|\nabla^2 v_{n,\e}(y)| > L\}}
\left|\e f\left(\theta_{n,\e}(y), \frac{1}{\e}\nabla^2 v_{n,\e}(y)\right)- 
f^{\infty}(\theta_{n,\e}(y),\nabla^2 v_{n,\e}(y))\right| \, dy \\
& \leq \int_{Q_{\nu}\cap \{\frac{1}{\e}|\nabla^2 v_{n,\e}(y)| > L\}}
C \, \e^{\gamma} \, |\nabla^2 v_{n,\e}(y)|^{1-\gamma} \, dy \\
& \leq C \, \e^{\gamma}
\left( \int_{Q_{\nu}}|\nabla^2 v_{n,\e}(y)| \, dy\right)^{1-\gamma}\\
& \leq C \, \e^{\gamma}\left(\int_{Q_{\nu}}
\e f\left(\theta_{n,\e}(y),\frac{1}{\e}\nabla^2 v_{n,\e}(y)\right) \, dy\right)^{1-\gamma} = O(\e^{\gamma}),
\end{align*}
since the integral in the last expression is uniformly bounded by \eqref{secondg}.
Hence, \eqref{limitffinfty} is proved. 

Therefore, from \eqref{secondg}, via a diagonalisation argument, we conclude that 
\begin{align*}
g(a,b,C,D,\nu) = \lim_{\varepsilon \to 0^+} \left[\displaystyle\int_{Q_{\nu}}f^\infty(\theta_{n(\e),\varepsilon}(y),
\nabla^2 v_{n(\e),\varepsilon}(y)) \, dy
+ |D\theta_{n(\e),\varepsilon}|(Q_{\nu})\right] 
\geq \widetilde K(a,b,C,D,\nu). 
\end{align*}

In order to prove the reverse inequality, let 
$\{\chi_n\} \subset BV(Q_\nu;\{0,1\})$, $\{u_n\} \subset W^{2,1}(Q_\nu;\Rb^d)$ be such that
$\chi_n \to \chi_{a,b,\nu}$ in $L^1(Q_\nu;\{0,1\})$, $u_n \to u_{C,D,\nu}$ in $W^{1,1}(Q_\nu;\Rb^d)$ and 
$${\widetilde K}(a,b,C,D,\nu) = 
\lim_{n\rightarrow+\infty}
\left[ \int_{Q_\nu} f^\infty(\chi_n(y),\nabla^2 u_n(y)) \, dy + |D\chi_n|(Q_\nu)\right].$$

For $x \in Q_\nu(0,\e)$, set 
$\displaystyle \theta_{n,\varepsilon}(x) := \chi_n\left(\frac{x }{\e}\right)$ and 
$\displaystyle v_{n,\varepsilon}(x) := \varepsilon u_n\left(\frac{x}{\e}\right)$. 
Then, changing variables and using the positive homogeneity of $f^\infty(q,\cdot)$, we obtain
\begin{align}\label{Kg1}
\widetilde{K}(a,b,C, D, \nu) & = \lim_{n\rightarrow+\infty}
\left[\int_{Q_\nu} f^\infty(\chi_n(y),\nabla^2 u_n(y)) \, dy + |D\chi_n|(Q_\nu)\right] \nonumber \\
& = \frac{1}{\e^{N-1}}\lim_{n\rightarrow+\infty}
\left[ \int_{Q_\nu(0,\e)} f^\infty(\theta_{n,\varepsilon}(x),
\nabla^2 v_{n,\varepsilon}(x)) \, dx + |D\theta_{n,\varepsilon}|(Q_\nu(x_0,\e))\right] \nonumber \\
& \geq \frac{1}{\e^{N-1}}\liminf_{n\rightarrow+\infty}
\left[\int_{Q_\nu(0,\e)} f(\theta_{n,\varepsilon}(x),\nabla^2 v_{n,\varepsilon}(x)) \, dx + |D\theta_{n,\varepsilon}|(Q_\nu(0,\e))\right] \nonumber \\
& \;\;\; + \frac{1}{\e^{N-1}}\liminf_{n\rightarrow+\infty}
\int_{Q_\nu(0,\e)} \Big(f^\infty(\theta_{n,\varepsilon}(x),
\nabla^2 v_{n,\varepsilon}(x)) 
- f(\theta_{n,\varepsilon}(x),\nabla^2 v_{n,\varepsilon}(x))\Big)\, dx =: I_1 + I_2.
\end{align}
Given that $\chi_n \to \chi_{a,b,\nu}$ in $L^1(Q_\nu;\{0,1\})$ and 
$u_n \to u_{C,D,\nu}$ in $W^{1,1}(Q_\nu;\Rb^d)$, it follows that, for every $\varepsilon$, $\theta_{n,\varepsilon} \to \chi_{a,b,\nu}$ in $L^1(Q_\nu(0,\e);\{0,1\})$ and $v_{n,\varepsilon} \to u_{C,D,\nu}$ in $W^{1,1}(Q_\nu(0,\e);\Rb^d)$, as $n\to +\infty$. Thus, 
\begin{equation}\label{Kg2}
	I_1 \geq \frac{1}{\e^{N-1}}
	\mathcal F(\chi_{a,b,\nu},u_{C,D,\nu};Q_\nu(0,\e)) 
	\geq \frac{1}{\e^{N-1}}
	m(\chi_{a,b,\nu},u_{C,D,\nu};Q_\nu(0,\e)). 
\end{equation}
On the other hand, the same calculations that were used to prove \eqref{limitffinfty} by means of hypothesis \eqref{finfty}, but using in this case the equivalent form given in \eqref{finfty2}, allow us to deduce that
\begin{equation}\label{Kg3}
\limsup_{\varepsilon \to 0^+}I_2 = 0.
\end{equation}
Hence, from \eqref{Kg1}, \eqref{Kg2} and \eqref{Kg3}, taking also into account
Remark \ref{remGM} ii), we conclude that
\begin{align*}
{\widetilde K}(a,b,C,D,\nu) 
\geq \limsup_{\varepsilon \to 0^+} \frac{1}{\e^{N-1}}
m(\chi_{a,b,\nu},u_{C,D,\nu};Q_\nu(0,\e)) 
 = g(a,b,C,D,\nu),
\end{align*}
and this completes the proof.
\end{proof}

We recall the notion of $BV$-ellipticity, given in \cite{AB}, and that of 
$SBH$-ellipticity defined in \cite{SZ2}.

\begin{Definition}\label{BVelldef}
Let $T \subset \mathbb R^d$ be a finite set. A function
$\Phi: T \times T \times \mathbb S^{N-1} \to [0,+\infty)$ is said to be
$BV$-elliptic if the following inequality holds
$$\Phi(a,b,\nu) \leq \int_{Q_\nu \cap S_\eta}\Phi(\eta^+,\eta^-,\nu_\eta) \, d \mathcal H^{N-1},$$
for every $a, b \in T$, $\nu \in \mathbb S^{N-1}$ and for every piecewise constant function $\eta \in SBV(Q_\nu;T)$ such that $\eta=a$ in 
$\partial Q_\nu \cap \{x \cdot \nu > 0\}$
and
$\eta=b$ in $\partial Q_\nu \cap \{x \cdot \nu < 0\}$.
\end{Definition}

Observe that when $T=\{a,b\}$, then Definition \ref{BVelldef} coincides with  $BV$-ellipticity with respect to $SJ_0$, as recently introduced in \cite[Definition 3.1]{EKM}. On the other hand, when dealing with perimeters, that is, when 
$a, b \in \{0,1\}$, among all the notions considered in
\cite[Definition 3.1]{EKM},  this is the only meaningful one. 

\begin{Remark}\label{Bvellch}
In view of \cite[Theorem 5.14]{AFP} and standard arguments in the Calculus of Variations which, in $BV$, allow us to prescribe the boundary conditions of the approximating sequences so as to coincide with those of the target function, 
$BV$-ellipticity can be characterised as follows.

A function
$\Phi: T \times T \times \mathbb S^{N-1} \to [0,+\infty)$ is 
$BV$-elliptic if
\begin{equation}\label{seqBVell}
\Phi(a,b,\nu) \leq \liminf_n\int_{Q_\nu \cap S_{\eta_n}}
\Phi(\eta_n^+,\eta_n^-,\nu_{\eta_n}) \, d \mathcal H^{N-1},
\end{equation}
for every $a, b \in T$, $\nu \in \mathbb S^{N-1}$ and for every sequence $\{\eta_n\} \subset SBV(Q_\nu; T)$ such that $\eta_n \to \eta_{a,b,\nu}$ in 
$L^{1}(Q_\nu;\mathbb R^d)$, $\eta_n = \eta_{a,b,\nu}$ on $\partial Q_\nu$, 
where $\eta_{a,b,\nu}$
is such that  $\eta=a$ in $\partial Q_\nu \cap \{x \cdot \nu > 0\}$
and $\eta=b$ in $\partial Q_\nu \cap \{x \cdot \nu < 0\}$.
\end{Remark}

\begin{Definition}\label{BHelldef}
Let $\mathcal P := \left \{(C,D,\nu) \in 
\mathbb R^{d \times N} \times \mathbb R^{d \times N} \times \mathbb S^{N-1} : 
C - D = P \otimes \nu, \mbox{ for some } P \in \mathbb R^d\right \}$.
A continuous function
$\varphi: \mathcal P \to [0,+\infty)$ is said to be
$SBH$-elliptic if the following inequality holds
\begin{equation}\label{seqSBHell}
\varphi(C,D,\nu) \leq \liminf_n\int_{Q_\nu \cap S_{\nabla u_n}}
\varphi((\nabla u_n)^+,(\nabla u_n)^-,\nu_{\nabla u_n}) \, d \mathcal H^{N-1},
\end{equation}
for every $(C,D,\nu) \in \mathcal P$ and for every sequence $\{u_n\} \subset SBH(Q_\nu;\mathbb R^d)$ such that $u_n \to u_{C,D,\nu}$ in 
$W^{1,1}(Q_\nu;\mathbb R^d)$, $u_n = u_{C,D,\nu}$ on $\partial Q_\nu$ and
$|\nabla^2 u_n| \to 0$ in $L^{1}(Q_\nu)$, \color{black} where $u_{C,D,\nu}$ is given by \eqref{bounddata}.
\end{Definition}

We end this section by proving a further property of our relaxed surface energy density $g$.

\begin{Proposition}\label{prop4.12}
Let $f$ be given by \eqref{densityint}, where $W_0$ and $W_1$ are continuous functions satisfying \eqref{growthint}, and assume, in addition, that $W_0$ and $W_1$ are positively homogeneous of degree one. 
If the function $g$ in \eqref{1gdef}, is continuous then it is 
is $BV\times SBH$-elliptic, in the sense that both \eqref{seqBVell} and \eqref{seqSBHell} hold,
that is, for every $a,b \in \{0,1\}$, $\nu \in \mathbb S^{N-1}$, $C,D \in \mathbb R^{d \times N}$ such that  
$C - D = P \otimes \nu$, for some $P \in \mathbb R^d$,
\begin{align*}
g(a,b,C,D,\nu) &=
\int_{Q_\nu}g(a,b, C, D, \nu) \, d \mathcal H^{N-1}(x) \\
&\leq \liminf_{n\to +\infty} \int_{Q_\nu} 
g(\chi^+_n(x),\chi^-_n(x), (\nabla u_n)^+(x),(\nabla u_n)^-(x), 
\nu_{\chi_n,\nabla u_n}(x)) \,d {\mathcal H}^{N-1}(x),
\end{align*}
whenever $\{\chi_n\} \subset BV(Q_\nu;\{0,1\})$, $\chi_n \to \chi_{a,b,\nu}$ in $L^1(Q_\nu;\{0,1\})$ and $\{u_n\} \subset SBH(Q_\nu;\mathbb R^d)$ is such that
$u_n \to u_{C,D,\nu}$ in $W^{1,1}(Q_\nu;\mathbb R^d)$, $|\nabla^2 u_n|\to 0$ in $L^1(Q_\nu)$ and $u_n= u_{C,D,\nu}$  on $\partial Q_\nu$, where $\chi_{a,b,\nu}$
and $u_{C,D,\nu}$ are given in \eqref{bounddata}. 
\end{Proposition}
\begin{proof}
Notice first that, due to the positive $1$-homogeneity of $W_0$ and $W_1$, we have
$$Q^2f(\chi,0)= \chi Q^2W_1(0)+ (1-\chi)Q^2W_0(0)\leq \chi W_1(0)+(1-\chi)W_0(0)=0,$$ 
so that $Q^2f(\chi,0) = 0$, for every $\chi \in BV(Q_\nu;\{0,1\})$.

Let $\{\chi_n\}$, $\{u_n\}$ be sequences satisfying the conditions of the statement.
Then, applying Theorem \ref{main} twice, first to $\chi_{a,b,\nu}$ and $u_{C,D,\nu}$ and then to $\chi_n$ and $u_n$, in $Q_\nu$, and by the lower semi-continuity of $\mathcal F$ (see Proposition \ref{firstprop}), it follows that
\begin{align*}
\int_{Q_\nu}g(a,b, C, D, \nu) \, d \mathcal H^{N-1}(x)
& =\int_{Q_\nu}Q^2 f(\chi_{a,b,\nu}(x), 0) \, dx + 
	\int_{Q_\nu}g(a,b, C, D, \nu) \, d \mathcal H^{N-1}(x)\\
&= \mathcal F(\chi_{a,b,\nu},u_{C,D,\nu}; Q_\nu)\\
&\leq \liminf_{n\to +\infty}\mathcal F(\chi_n,u_n; Q_\nu)\\
&\leq \liminf_{n\to +\infty}\left(\int_{Q_\nu} 
Q^2f(\chi_n(x),\nabla^2 u_n(x)) \, dx \right.\\
& \left.\hspace{1cm} + \int_{Q_\nu} 
g(\chi^+_n(x), \chi^-_n(x), (\nabla u_n)^+(x),(\nabla u_n)^-(x), \nu_{\chi_n,\nabla u_n}(x)) \, d {\mathcal H}^{N-1}(x)\right)\\
&= \liminf_{n\to +\infty} \int_{Q_\nu} 
g(\chi^+_n(x), \chi^-_n(x), (\nabla u_n)^+(x),(\nabla u_n)^-(x), \nu_{\chi_n,\nabla u_n}(x)) \, d {\mathcal H}^{N-1}(x),
\end{align*}
where we used Lebesgue's dominated convergence theorem in the last equality.

\end{proof}

\appendix

\section{}\label{app}

In this appendix we provide another partial characterisation of the relaxed surface energy density $g$. To this end, we introduce an auxiliary functional,
similar to \eqref{mbis}, but considering more regular functions. 

For $\chi \in BV(\O;\{0,1\})$, $u \in BH(\O;\mathbb R^d)$ and 
$A \in \mathcal O_{\infty}(\O)$  let 
\begin{align}\label{m0}
		m_0(\chi,u;A) := &\inf\Big\{F (\theta,v;A) : 
\theta \in  BV(A;\{0,1\}), v \in W^{2,1}(A;\Rb^d), 
\theta \equiv \chi  \text{ on } \partial A,
v \equiv u \text{ on } \partial A \Big\}.
\end{align}

The following lemma shows that $m$ is bounded from below by $m_0$.

\begin{Lemma}\label{BFT314}
Let $f$ be given by \eqref{densityint}, where $W_0$ and $W_1$ are continuous functions satisfying \eqref{growthint}.
Let $m$ and $m_0$ be given by \eqref{mbis} and \eqref{m0}, respectively. Then, for every $\chi \in BV(\O;\{0,1\})$, $u \in BH(\O;\mathbb R^d)$ and 
$A \in \mathcal O_{\infty}(\O)$, we have
$m(\chi,u;A) \geq m_0(\chi,u;A).$
\end{Lemma}
\begin{proof}
To show the desired inequality, we fix $\e > 0$ and let 
$\theta \in  BV(\O;\{0,1\}), v \in BH(\O;\mathbb R^d)$ be such that 
$\theta \equiv \chi$ on $\partial A$, ${\rm supp}(v - u) \Subset A$ and
$$m(\chi,u;A) + \e \geq \mathcal F(\theta,v;A).$$
By Proposition~\ref{firstprop} ii), let $\{\chi_n\} \subset  BV(A;\{0,1\}), \{u_n\} \subset W^{2,1}(A;\Rb^d)$ satisfy $\chi_n \to \theta$ in $L^1(A;\{0,1\})$, 
$u_n \to v$ in $W^{1,1}(A;\Rb^d)$ and
$$\mathcal F(\theta,v;A) = \lim_{n \to + \infty}F(\chi_n,u_n;A).$$
Lemma~\ref{Lemma1KRBH} ensures the existence of a regular sequence
$\{v_n\} \subset W^{2,1}(A;\Rb^d)$ such that $v_n \to v$ in $W^{1,1}(A;\Rb^d)$,
$v_n = v = u \mbox{ on } \partial A$ and $\nabla^2v_n \to D^2v$ (area) $\langle \cdot \rangle$-strictly in $A$.
We now apply Proposition~\ref{newslicing} and Remark \ref{remreg}, to conclude that there exists a subsequence
$\{v_{n_k}\}$ of $\{v_n\}$ and there exist sequences 
$\{w_k\} \subset W^{2,1}(A;\Rb^d)$,
$\{\eta_k\} \subset BV(A;\{0,1\})$ verifying $ {\rm supp}(w_k -v_{n_k})  \Subset A$,
$\eta_k \equiv \theta \equiv \chi$ on  $\partial A$ and
$$\limsup_{k \to +\infty}F(\eta_k,w_k;A) \leq 
\liminf_{n \to +\infty}F(\chi_n,u_n;A).$$
Therefore, the pair $(\eta_k,w_k)$ is admissible for $m_0(\chi,u;A)$ so that
$$m_0(\chi,u;A) \leq \limsup_{k \to +\infty}F(\eta_k,w_k;A) \leq 
\liminf_{n \to +\infty}F(\chi_n,u_n;A) = \mathcal F(\theta,v;A) \leq 
m(\chi,u;A) + \e,$$
so the desired inequality follows by letting $\e \to 0^+.$
\end{proof}

For $(a,b,C,D,\nu)  \in \{0,1\} \times \{0,1\} \times \mathbb{R}^{d\times N} \times\mathbb{R}^{d\times N} \times \mathbb S^{N-1},$ such that 
$C-D= P \otimes \nu$ for some $P \in \mathbb R^d$, define
\begin{equation}\label{KSQ}
 K(a,b,C,D,\nu):=\inf\left\{  
 	\displaystyle\int_{Q_{\nu}}\!\!\!
 	(Q^2f)^{\infty}(\chi(x), \nabla^2 u(x)) \, dx+|D\chi|(Q_{\nu}):\left(\chi,u\right)
 	\in\mathcal{A}(a,b,C,D,\nu)\right\}, 
 \end{equation}
where the set of admissible functions is given by
\begin{align}\label{admA}
\mathcal{A}(a,b,C,D,\nu) :=\bigg\{ \left(\chi,u\right)\in
BV\left(Q_{\nu};\{0,1\}\right)  \times 
W^{2,1}\left(Q_{\nu};\mathbb{R}^{d}\right): 
\chi= \chi_{a,b,\nu} \text{ on } \partial Q_\nu, 
u = u_{C,D,\nu} \text{ on }\partial Q_\nu\bigg\}, 
\end{align}
and $\chi_{a,b,\nu}$, $u_{C,D,\nu}$ were defined in \eqref{bounddata}.

Taking into account Remark \ref{remGM} ii), in the following result we show that $K$ is a lower bound for $g$. 
The assumption that $f$ is $2$-quasiconvex in the second variable is placed merely for simplicity since, in fact, we prove that the inequality $g\geq K$ holds if, in \eqref{KSQ}, $(Q^2 f)^\infty$  is replaced by $f^\infty$.

\begin{Proposition}\label{prop4.9}
Let $W_0$ and $W_1$ be as in \eqref{growthint} and let $f$ be as in \eqref{densityint}.
Assume that $f$ is $2$-quasiconvex in the second variable and that 
\eqref{finfty} holds, where $f^\infty$ is defined in \eqref{recS}.
Given $\chi \in BV(\Omega;\{0,1\})$ and $u \in BH(\Omega;\mathbb R^d)$,  
for $\mathcal H^{N-1}$ a.e. $x_0 \in \Omega \cap (S_\chi \cup S_{\nabla u})$, we have
$$g(\chi^+(x_0),\chi^-(x_0),(\nabla u)^+(x_0),(\nabla u)^-(x_0),\nu(x_0)) \geq
K(\chi^+(x_0),\chi^-(x_0),u^+(x_0),u^-(x_0),\nu(x_0)),$$
where $\chi^+(x_0) = \chi^-(x_0) = \widetilde\chi(x_0)$ if 
$x_0 \in S_{\nabla u} \setminus S\chi$ and 
$(\nabla u)^+(x_0) = (\nabla u)^-(x_0) = \widetilde{\nabla u}(x_0)$ if
$x_0 \in S_\chi \setminus S_{\nabla u}$.
\end{Proposition}
\begin{proof}
Let  $x_0 \in \Omega \cap (S_\chi \cup S_{\nabla u})$. As before, for simplicity of notation, we write $\nu = \nu(x_0)$.
	
By Proposition \ref{mg} and Lemma \ref{BFT314}, we have
\begin{align*}
&g(\chi^+(x_0),\chi^-(x_0),(\nabla u)^+(x_0), (\nabla u)^-(x_0),\nu) \\
& \qquad \geq\limsup_{\e \to 0^+}\frac{1}{\e^{N-1}}
\inf\Big\{F (\theta,v;Q_\nu(0,\e)) : 
\theta \in  BV(Q_\nu(0,\e);\{0,1\}), v \in W^{2,1}(Q_\nu(0,\e);\Rb^d), \\
& \hspace{4,5cm} \theta = \chi_0
\mbox{ on } \partial Q_\nu(0,\e), v = u_0 
\mbox{ on } \partial Q_\nu(0,\e)\Big\},
\end{align*}
where, for $x \in Q_\nu(0,\varepsilon)$, $\chi_0$ and $u_0$ are given by
\begin{equation*}
\chi_0(x) := \begin{cases}
\chi^+(x_0), & {\rm if } \; x \cdot \nu > 0,\\
\chi^-(x_0), & {\rm if } \; x \cdot \nu < 0
\end{cases} \quad \mbox{and} \quad
u_0(x) := \begin{cases}
(\nabla u)^+(x_0)\cdot x, & {\rm if } \; x \cdot \nu > 0,\\
(\nabla u)^-(x_0)\cdot x, & {\rm if } \; x \cdot \nu < 0.
\end{cases}
\end{equation*}
Therefore, for every $n \in \mathbb N$, there exist 
$\{\theta_{n,\e}\} \subset  BV(Q_\nu(0,\e);\{0,1\})$, 
$\{v_{n,\e}\} \subset W^{2,1}(Q_\nu(0,\e);\Rb^d)$ such that $\theta_{n,\e} = \chi_0$ on $\partial Q_\nu(0,\e)$, $v_{n,\e} = u_0$  on $\partial Q_\nu(0,\e)$ and 
\begin{align}\label{g}
& g(\chi^+(x_0),\chi^-(x_0),(\nabla u)^+(x_0),(\nabla u)^-(x_0),\nu) + \frac{1}{n} \nonumber\\
&\hspace{2cm}\geq \limsup_{\e \to 0^+}\frac{1}{\e^{N-1}}
\left[\int_{Q_{\nu}(0,\e)}f(\theta_{n,\e}(x), \nabla^2 v_{n,\e}(x)) \, dx + |D\theta_{n,\e}|(Q_{\nu}(0,\e))\right] \nonumber\\
& \hspace{2cm}= \limsup_{\varepsilon \rightarrow 0^+}
\left[\int_{Q_{\nu}}\e f(\theta_{n,\e}(\e y), \nabla^2 v_{n,\e}(\e y)) \, dy + 
\int_{Q_{\nu}\cap \frac{1}{\e}S_{\theta_{n,\e}}}d\mathcal H^{N-1}(y)\right] \nonumber\\
& \hspace{2cm}= \limsup_{\varepsilon \rightarrow 0^+}
\left[\int_{Q_{\nu}}\e f\left(\chi_{n,\e}(y), \frac{1}{\e}\nabla^2 u_{n,\e}(y)\right) \, dy + 
|D\chi_{n,\e}|(Q_{\nu})\right]  \nonumber\\
& \hspace{2cm} \geq  \liminf_{\varepsilon \rightarrow 0^+}
\left[\int_{Q_{\nu}}f^{\infty}(\chi_{n,\e}(y),\nabla^2 u_{n,\e}(y)) \, dy + 
|D\chi_{n,\e}|(Q_{\nu})\right]  \nonumber\\
& \hspace{2cm} +  \liminf_{\varepsilon \rightarrow 0^+}
\int_{Q_{\nu}}\left[\e f\left(\chi_{n,\e}(y), \frac{1}{\e}\nabla^2 u_{n,\e}(y)\right)  -
f^{\infty}(\chi_{n,\e}(y),\nabla^2 u_{n,\e}(y))\right]  dy,
\end{align}
where, for $y \in Q_\nu$, $\chi_{n,\e}(y) = \theta_{n,\e}(\e y)$ and 
$\displaystyle u_{n,\e}(y) = \frac{1}{\varepsilon}v_{n,\e}(\e y)$. 
Arguing as in the proof of \eqref{limitffinfty}, we conclude that 
\begin{equation}\label{limit0}
\liminf_{\e \to 0^+}
\int_{Q_{\nu}}\left[\e f\left(\chi_{n,\e}(y), \frac{1}{\e}\nabla^2 u_{n,\e}(y)\right)- f^{\infty}(\chi_{n,\e}(y),\nabla^2 u_{n,\e}(y))\right]  dy  = 0.
\end{equation}
Thus, since $\{(\chi_{n,\e},u_{n,\e})\} \subset \mathcal{A}(\chi^+(x_0),\chi^-(x_0),(\nabla u)^+(x_0),(\nabla u)^-(x_0),\nu)$, where this class of admissible fields is given in \eqref{admA}, 
from \eqref{g}, \eqref{limit0} and the definition of $K(\chi^+(x_0),\chi^-(x_0),(\nabla u)^+(x_0),(\nabla u)^-(x_0),\nu)$ in \eqref{KSQ},
it follows that
\begin{align*}
&g(\chi^+(x_0),\chi^-(x_0),(\nabla u)^+(x_0),(\nabla u)^-(x_0),\nu) + \frac{1}{n} \\
&\hspace{3cm}
\geq \liminf_{\e \to 0^+}
\left[\int_{Q_{\nu}}f^{\infty}(\chi_{n,\e}(y),\nabla^2 u_{n,\e}(y)) \, dy + 
|D\chi_{n,\e}(Q_{\nu})\right] \\
&\hspace{3cm}\geq K(\chi^+(x_0),\chi^-(x_0),u^+(x_0),u^-(x_0),\nu),
\end{align*}
hence the result follows by letting $n \to +\infty$.
\end{proof}

\section*{Acknowledgements}
The research of ACB was partially supported by National Funding from FCT -
Funda\c c\~ao para a Ci\^encia e a Tecnologia through project 
UIDB/04561/2020: https://doi.org/10.54499/UIDB/04561/2020.
She also acknowledges the support of GNAMPA-INdAM through Programma Professori Visitatori 2023.
EZ is a member of GNAMPA-INdAM, whose support is gratefully acknowledged through project GNAMPA-INdAM 2023 ``Prospettive nelle scienze dei materiali: modelli variazionali, analisi asintotica e omogeneizzazione'' CUP E53C22001930001 and PRIN ``Mathematical Modeling of Heterogeneous Systems'', CUP 853D23009360006. 
She is also very grateful to CMAFcIO at the University of Lisbon for its kind hospitality and support.

The authors thank Prof. Martin Kru\v zík for his valuable comments, which improved this manuscript.

\end{document}